\documentclass{article}

\usepackage{amsmath,amssymb,amsthm}
\usepackage{latexsym}
\usepackage[T1]{fontenc}
\usepackage[francais]{babel}

\def \vs{\vskip .6cm}

\def \beq{\begin{eqnarray*}}
\def \eeq{\end{eqnarray*}}

\def \n{\nabla}
\def \nt{\widetilde{\nabla}}
\def \nb{\overline \nabla}
\def \Rt{\widetilde R}
\def \Rb{\overline R}

\def \RM{\mathbb{R}}

\def \CM{\mathbb{C}}

\def \J {\mathcal J}

\newcommand \bo [1] {\bigotimes \phantom{ }^#1}

\def \LLambda {{\bf \Lambda}}
\def \llambda {\boldsymbol{\lambda}}

\def \leftr {[\hbox{\hspace{-0.15em}}[}
\def \rightr {]\hbox{\hspace{-0.15em}}]}
\def \la {\langle}
\def \ra {\rangle}

\newtheorem{defi}{D\'efinition}[section]

\newtheorem{propdefi}[defi]{Proposition et d\'efinition}

\newtheorem{prop}[defi]{Proposition}

\newtheorem{theo}[defi]{Th\'eor\`eme}

\newtheorem{theo'}{Th\'eor\`eme}

\newtheorem{lemm}[defi]{Lemme}

\newtheorem*{NB}{Remarque}

\newtheorem{coro}[defi]{Corollaire}

\title{Espace de twisteurs d'une vari\'et\'e presque
  hermitienne de dimension 6}

\date{}

\author{Jean-Baptiste Butruille}

\begin{document}

\maketitle

\begin{abstract}

On s'int\'eresse \`a l'espace de twisteurs
r\'eduit d'une vari\'et\'e presque hermitienne, en relisant un article
de N.R. O'Brian et J.H. Rawnsley \cite{ob}. On traite la question laiss\'ee
ouverte de la dimension 6. Cet espace est muni d'une
structure presque complexe $\J$ en utilisant la distribution
horizontale de la connexion hermitienne canonique. On montre qu'une condition n\'ec\'essaire d'int\'egrabilit\'e de $\J$ est que la vari\'et\'e soit de type $W_1
\oplus W_4$ dans la classification de Gray et Hervella \cite{gr4}. Dans la
deuxi\`eme partie on montre alors que les seules vari\'et\'es de type
$W_1 \oplus W_4$ en dimension 6 sont les vari\'et\'es localement
conform\'ement \og nearly K\"ahler \fg. Finalement la
structure presque complexe de l'espace de twisteurs r\'eduit est
int\'egrable si et seulement si la
vari\'et\'e est localement conforme \`a la sph\`ere $S^6$ ou \`a
une vari\'et\'e k\"ahlerienne, Bochner-plate.
\end{abstract}

\section*{Introduction}

La th\'eorie des twisteurs invent\'ee par R. Penrose (voir l'article
fondateur \cite{pe}) est un moyen d'utiliser les techniques efficaces
de la g\'eom\'etrie holomorphe pour r\'esoudre des probl\`emes de
g\'eom\'etrie riemannienne ou pseudo-riemannienne.

Soit $M$ une vari\'et\'e de dimension paire $m=2n$. On part d'une
vari\'et\'e complexe $Z$, donn\'ee avec une submersion \`a fibres
complexes $\pi : Z \to M$. On associe \`a tout point $j$ de $Z$ un
endomorphisme de carr\'e -1 de $T_{\pi(j)}M$ (ou \`a toute
section, une structure presque complexe de $M$) en transportant la
multiplication par $i$ de $T_j Z$ par l'isomorphisme d\'ependant
du point $(\pi_*)_j : T_j Z / \mathcal V_j \to T_{\pi(j)}M$, o\`u $\mathcal V$ est
la distribution verticale. Soit $\mathcal Z$ le fibr\'e de $M$
dont la fibre au-dessus de $x$ est l'ensemble des endomorphismes
de carr\'e -1 de $T_x M$. On note $\pi_0 : \mathcal Z \to M$ la
projection canonique. On a donc une application $\varphi$ de $Z$
dans $\mathcal Z$, pr\'eservant les fibres. D'autre part, il est
connu que $\mathcal Z_x = \pi_0^{-1}(x)$ est isomorphe en tout
point \`a l'espace sym\'etrique hermitien $GL(m,\RM)/GL(n,\CM)$
et admet par cons\'equent une structure presque complexe naturelle int\'egrable.
Alors on demande que $\varphi$ soit injective et que pour tout $x
\in M$ la restriction de $\varphi$ \`a $Z_x$ soit une injection
holomorphe. Dans ce cas, $Z$ est appel\'e un {\it espace de twisteurs
complexe} de $M$.

R\'eciproquement, pour obtenir un espace de twisteurs complexe sur
$M$, on prend $Z$ une sous-vari\'et\'e de $\mathcal Z$ telle que la restriction de
$\pi_0$ \`a $Z$ est toujours une fibration et pour tout $x \in M$,
$Z_x$ est une sous-vari\'et\'e complexe de $\mathcal Z_x$. On construit une
structure presque complexe $\J$ sur $Z$, en se servant d'une section de
la suite exacte
\[ 0 \to V \to TZ \to TZ/\mathcal V \to 0
\]
donn\'ee
d'habitude par une connexion sur $M$. Alors $Z$ est un espace de
twisteurs complexe, avec $\varphi$ l'injection
canonique, si et seulement si $\J$ est int\'egrable.

Le parfait exemple d'une telle situation est la fibration \`a fibres
$\CM P^1$ de l'espace projectif complexe $\CM P^3$ sur la
sph\`ere $S^4$ dont les vertus furent
d\'ecouvertes par Atiyah, Hitchin et Singer dans \cite{at}. Les m\^emes
auteurs ont cherch\'e une g\'en\'eralisation aux vari\'et\'es
riemanniennes de dimension 4, en posant a priori que $Z$ est la
sous-vari\'et\'e de $\mathcal Z$ constitu\'ee des structures
presque complexes compatibles avec la m\'etrique. Cela peut d'ailleurs
\^etre fait en
dimension paire sup\'erieure mais d\`es la dimension 6 la
condition obtenue pour l'int\'egrabilit\'e de $\J$ est que la
vari\'et\'e soit conform\'ement plate tandis qu'en dimension 4
elle a lieu pour toute la riche classe des vari\'et\'es
auto-duales, en raison d'une singularit\'e de la d\'ecomposition
en composantes irr\'eductibles, en $m=4$, de la r\'epr\'esentation de $SO(m)$ sur
l'espace des tenseurs de courbure riemannienne abstraits.

O'Brian et Rawnsley \cite{ob} regardent, eux, des espaces de twisteurs
associ\'es \`a une $G$-structure et une $G$-connexion. Le cas
originel correspond bien s\^ur \`a $G=SO(m)$ et la connexion de
Levi-Civita. Ils se sont particuli\`erement int\'eress\'es au cas
o\`u $G=U(n)$. A leur suite, on consid\`ere une vari\'et\'e presque
hermitienne $(M,g,J_0)$. On demande que les sections de $Z$ soient
compatibles avec
$g$ et commutent avec $J_0$ et on a besoin pour construire $\J$ d'une
connexion hermitienne $\nt$. Les
conditions d'int\'egrabilit\'e rappel\'ees section 2 portent alors non
seulement sur la courbure de $\nt$ mais sur sa
torsion, qui n'est pas nulle en g\'en\'eral. Par ailleurs on montre
section 3 qu'on peut sans perte de g\'en\'eralit\'e pour notre
probl\`eme choisir la connexion hermitienne canonique $\nb$.

En dimension sup\'erieure \`a $10$, $\J$ est int\'egrable si et
seulement si la vari\'et\'e est localement conforme \`a une
vari\'et\'e k\"ahlerienne (LCK) dont le tenseur de Bochner est
nul. Cette classe de vari\'et\'es pr\'esente elle-m\^eme un grand int\'er\^et. Leur \'etude
difficile est abord\'ee par exemple dans \cite{br2} (voir aussi
\cite{da}). On s'int\'eresse
ici \`a la dimension 6. On montre que les conditions impos\'ees
\`a la torsion sont moins strictes en cette dimension puisque outre les vari\'et\'es
LCK, toutes les vari\'et\'es de type $W_1 \oplus W_4$ dans la
classification de Gray-Hervella \cite{gr4} les satisfont. Parmi celles-ci
on trouve en particulier les vari\'et\'es strictement \og nearly K\"ahler \fg \ 
(NK). Mais la forme de la courbure de ces derni\`eres est si particuli\`ere en dimension 6 que les conditions
impos\'ees \`a celle-ci (proposition \ref{O'Brian & Rawnsley}) ne laissent finalement que la
sph\`ere $S^6$.

Or on montre, section 5, que
\begin{theo'}
Les vari\'et\'es presque hermitiennes de type $W_1
\oplus W_4$ en dimension 6 sont localement conformes \`a une vari\'et\'e NK.
\label{W1+W4=W1}
\end{theo'}

Ce th\'eor\`eme et la discussion pr\'ec\'edente permettent de conclure, compte-tenu de l'invariance conforme de l'espace de twisteurs r\'eduit et de sa structure presque complexe :
\begin{theo'}
Soient $M$ une vari\'et\'e presque hermitienne de dimension 6, $Z$
son espace de twisteurs r\'eduit, $\J$ la structure presque
complexe sur $Z$ associ\'ee \`a la connexion hermitienne canonique
: $\J$ est int\'egrable si et seulement si $M$ est localement
conforme \`a une vari\'et\'e k\"ahlerienne Bochner-plate ou \`a la sph\`ere
$S^6$ munie de sa structure NK. \label{twisteurs}
\end{theo'}

L'article est organis\'e en deux parties,
correspondant aux deux th\'eor\`emes principaux \ref{W1+W4=W1} et
\ref{twisteurs}, quasiment ind\'ependantes.
La r\'esolution compl\`ete du probl\`eme soulev\'e dans la premi\`ere partie a motiv\'e l'\'ecriture de la seconde partie.

La m\'ethode utilis\'ee \`a la section 5 s'inspire de l'\'etude des
vari\'et\'es NK de dimension 6. A. Gray a montr\'e dans \cite{gr3} que
celles-ci sont soit k\"ahleriennes, soit {\it strictement} NK (SNK). Dans
le dernier cas elles admettent une r\'eduction naturelle \`a $SU(3)$. C'est de
tenir toujours un meilleur compte de cette structure $SU(3)$ que
sont venus les derniers r\'esultats les concernant. D'abord
Reyes-Carrion \cite{re} a montr\'e que la connexion hermitienne
canonique \'etait
en fait une connexion $SU(3)$. Puis il a d\'ecouvert, ce que
Hitchin a rendu explicite dans \cite{hi3}, que toute l'information pour la
structure $SU(3)$, y compris la m\'etrique et la structure presque
complexe, est comprise dans la donn\'ee de deux formes : la forme de
K\"ahler $\omega$ et la forme volume complexe $\Psi$ (ou dans ce cas
la diff\'erentielle de la forme de K\"ahler $d\omega$), ce qui permet
de caract\'eriser les vari\'et\'es SNK en dimension 6 par une \'equation
diff\'erentielle simple portant sur la structure $SU(3)$.

Ici on s'int\'eresse \`a d'autres vari\'et\'es presque hermitiennes de
dimension 6 qu'on appelle {\it sp\'eciales} c'est-\`a-dire \`a d'autres
structures $U(3)$ qui induisent une structure $SU(3)$ sur la
vari\'et\'e par l'interm\'ediaire de $d\omega$. Salamon, Chiossi \cite{ch},
prolongeant le travail de Gray, Hervella ont classifi\'e les
vari\'et\'es $SU(3)$ en consid\'erant la torsion
intrins\'eque. Celle-ci est donn\'ee par la d\'ecomposition en types
de $d\omega$, $d\Psi$. Les vari\'et\'es $W_1 \oplus W_4$ sont alors
caract\'eris\'ees par deux \'equations diff\'erentielles portant notamment
sur la forme de Lee $\theta$ et on peut montrer que celle-ci est
ferm\'ee, c'est-\`a-dire repr\'esente localement (par le lemme de Poincar\'e)
un changement conforme de m\'etrique par lequel la vari\'et\'e est issue d'une
vari\'et\'e NK.

Section 6, vue l'invariance conforme de la d\'efinition de
l'espace de twisteurs, r\'eduit ou non, on reformule les r\'esultats de la
section 5 en
consid\'erant des vari\'et\'es presque hermitiennes {\it conformes}. On
laisse ouverte la question de savoir si un th\'eor\`eme tel que
\ref{W1+W4=W1} a lieu en toute dimension et pour d'autres classes
de vari\'et\'es presque hermitiennes, stables par transformation
conforme. Cette question est li\'ee \`a l'existence des vari\'et\'es de
type $G_1$, $G_2$ de Hervella, Vidal \cite{he}. On donne seulement, section 7, un
r\'esultat d'existence locale de vari\'et\'es de type $W_1 \oplus W_2
\oplus W_4$ de dimension 6 non localement conformes \`a des vari\'et\'es
de type $W_1 \oplus W_2$.

\section{Pr\'eliminaires}

On souhaite donner ici quelques d\'efinitions g\'en\'erales et quelques r\'esultats simples ou classiques sur les vari\'et\'es presque hermitiennes.

D'abord fixons quelques notations. Soient $M$ une vari\'et\'e de dimension
$m$ et $GL(M)$ le fibr\'e principal sur $M$ de groupe $GL(m)$ (le fibr\'e
des rep\`eres). Les repr\'esentations de $GL(m)$ fournissent des fibr\'es
associ\'es de $GL(M)$, les fibr\'es de tenseurs de $M$, parmi lesquels :
$TM$, le fibr\'e tangent ou les fibr\'es ext\'erieurs $\LLambda^p$.

Maintenant, si $M$ est une vari\'et\'e riemannienne orient\'ee, on a une
premi\`ere r\'eduction de $GL(M)$ \`a $SO(m)$ : soit $SO(M)$ le fibr\'e des
rep\`eres orthonorm\'es directs de $M$. On note $\mathfrak{so}(M)$ son
fibr\'e adjoint c'est-\`a-dire le fibr\'e des endomorphismes antisym\'etriques
de $TM$. Plus g\'en\'eralement, soient $G$ un groupe de Lie, $G \subset
SO(m)$, d'alg\`ebre de Lie $\mathfrak g$ et $\mathfrak g^{\perp}$
l'orthogonal de $\mathfrak g$ dans $\mathfrak{so}(m)$.
On suppose que $M$
admet une r\'eduction \`a $G$, c'est-\`a-dire qu'il existe un sous-fibr\'e principal
$G(M)$ de groupe $G$ de $SO(M)$. Alors on
note $\mathfrak g(M)$ le fibr\'e adjoint de
$G(M)$. De m\^eme $\mathfrak g^{\perp}(M)$ d\'esignera le fibr\'e
associ\'e \`a la repr\'esentation $\mathfrak g^{\perp}$ de $G$.

Dans cet article on consid\`ere des vari\'et\'es presque hermitiennes. Une telle vari\'et\'e
est d\'efinie en dimension paire $m=2n$ par une r\'eduction du fibr\'e des
rep\`eres \`a $U(n)$, not\'ee $U(M)$, ou autrement par une m\'etrique $g$ et
une structure presque complexe $J$, orthogonale,
\[ \forall X,Y \in TM, \quad g(JX,JY)=g(X,Y),
\]
d\'efinissant une 2-forme $\omega$, appel\'ee \emph{forme de K\"ahler} :
\[ \forall X,Y \in TM, \quad \omega(X,Y)=g(JX,Y) \]
Soit $T^{1,0} \subset T^{\CM} M$ le sous-fibr\'e des vecteurs complexes de type (1,0) par rapport \`a $J$. De m\^eme, soit $T^{0,1}$ le fibr\'e des vecteurs de type (0,1). Autrement dit, quel que soit $x \in M$, $T^{1,0}_x$ (resp. $T^{0,1}_x$) est le sous-espace propre (complexe) de $J_x$ pour la valeur propre $i$ (resp. $-i$). Le tenseur de Nijenhuis $N$ mesure l'int\'egrabilit\'e de la structure presque complexe $J$ ou de la distribution $T^{1,0}$. Il est d\'efini, en vertu du th\'eor\`eme de Frobenius, par
\begin{equation}
\forall X,Y \in TM, \quad N(X,Y) + iJN(X,Y) = [X^{1,0},Y^{1,0}]^{0,1},
\label{N}
\end{equation}
o\`u pour tout $X \in T^{\CM} M$, $X^{1,0}=\frac{1}{2}(X-iJX)$ et $X^{0,1}=\frac{1}{2}(X+iJX)$ d\'esignent les projections de $X$ sur $T^{1,0}$, $T^{0,1}$, respectivement. Soit donc $(M,g,J)$ une vari\'et\'e presque hermitienne. Le fibr\'e $\mathfrak u(M)$ (resp. $\mathfrak u(M)^{\perp}$) est le
fibr\'e des endomorphismes antisym\'etriques de $TM$ qui commutent
(resp. anticommutent) \`a $J$.

On rappelle aussi que pour tout fibr\'e principal $P$ il existe une action
naturelle du fibr\'e int\'erieur sur les fibr\'es associ\'es. En particulier --
si $P$ est un sous-fibr\'e de $GL(M)$ -- sur
les fibr\'es de tenseurs. Il existe aussi une action d\'eriv\'ee du fibr\'e adjoint. On note $\Phi s$ (resp. $A.s$) l'action d'un automorphisme vertical $\Phi$ de $P$ (resp. d'une section $A$ du fibr\'e adjoint) sur un champ de tenseurs $s$. On
souhaite expliciter cette action dans le cas o\`u $s$ est un tenseur de type (2,1) : quels que soient $X,Y,Z \in TM$,
\[ \Phi s(X,Y) = \Phi\big(s(\Phi^{-1}X,\Phi^{-1}Y)\big) \]
\[ A.s(X,Y) = A s(X,Y) -s(AX,Y) -s(X,AY) \]
Or la structure presque complexe $J$ peut-\^etre vue alternativement
comme un endomorphisme orthogonal de $TM$ ou comme un endomorphisme
antisym\'etrique, c'est-\`a-dire une section du fibr\'e adjoint
$\mathfrak{so}(M)$. Par cons\'equent, elle peut agir de ces deux fa\c cons et puisque $J^2 = -Id$,
\[ J s(X,Y,Z) = J\big(s(JX,JY)\big) \]

Les repr\'esentations d'un groupe, ici $U(n)$, interviennent dans la
construction des fibr\'es associ\'es. On
utilise la notation de Salamon \cite{sa}. Soit $V$ un espace de
repr\'esentation
complexe. On dit que $V$ est de \emph{type complexe} si $V$ n'est
pas isomorphe \`a son conjugu\'e $\overline V$. Alors $\leftr V \rightr$
d\'esigne
simplement l'espace vectoriel (ou l'espace de repr\'esentation)
r\'eel sous-jacent. Au contraire $V$ est dit de \emph{type
r\'eel} s'il admet une structure
r\'eelle, soit un endomorphisme de carr\'e 1 anticommutant
\`a $i$. En particulier $V \simeq \overline V$. Alors on note
$[V]$ le sous-espace propre de cet endomorphisme pour la valeur
$1$. Il s'agit d'un espace vectoriel r\'eel qui admet $V$ comme
complexification : $V \simeq
[V] \otimes_{\RM} \CM$. Selon les cas on a donc $\text{dim}_{\RM}
\leftr V \rightr = 2 \text{dim}_{\CM} V$ ou $\text{dim}_{\RM} [V]
= \text{dim}_{\CM} V$. De plus on note $\leftr \bf V \rightr$ ou $[\bf V]$ le fibr\'e associ\'e de $U(M)$ correspondant.

L'exemple fondamental pour nous est l'espace des formes de type
$(p,q)$. On sait que l'intersection de cet espace $\lambda^{p,q}$ avec
les r-formes r\'eelles ($r=p+q$) est nulle si $p \neq q$. En
revanche $\lambda^{p,q} \oplus \lambda^{q,p}$ est la complexification
  de $\leftr \lambda^{p,q}
  \rightr$, l'espace des formes \emph{r\'eelles} de type $(p,q)+(q,p)$.
De l'autre c\^ot\'e si $p=q$ on a directement que $\lambda^{p,p}$ est
la complexification de $[\lambda^{p,p}]$.

Pour une vari\'et\'e riemannienne, la m\'etrique fournit un isomorphisme
$SO(m)$-invariant des deux fibr\'es $\mathfrak{so}(M)$ et
$\LLambda^2$. Pour une vari\'et\'e presque hermitienne $(M,g,J)$ on a en
outre les isomorphismes $U(n)$-invariants :
\[ \mathfrak{u}(M) \simeq [ \llambda^{1,1} ] \quad \text{et} \quad
\mathfrak{u}(M)^{\perp} \simeq \leftr \llambda^{2,0} \rightr \]

Soit $\n$ la connexion de Levi-Civita de $g$. Gray et Hervella
\cite{gr4} regardent la d\'eriv\'ee covariante de la forme de K\"ahler $\n
\omega$. C'est une section de $\LLambda^1 \otimes \leftr \llambda^{2,0} \rightr$. En effet $J^2=-Id$ implique
\begin{equation}
(\n_X J)J + J(\n_X J) = 0
\label{nJ}
\end{equation}
donc $\n J$ est une section de $\LLambda^1 \otimes
\mathfrak{u}(M)^{\perp}$. On peut voir ce tenseur comme le d\'efaut pour $M$ d'\^etre
k\"ahlerienne. En vue de classifier les vari\'et\'es
presque hermitiennes Gray et Hervella d\'ecomposent $\Lambda^1 \otimes
\leftr \lambda^{2,0} \rightr$ sous l'action de $U(n)$ :
\beq
\Lambda^1 \otimes \leftr \lambda^{2,0} \rightr & = & \leftr
\lambda^{1,0} \otimes \lambda^{2,0} \rightr \oplus \leftr
\lambda^{0,1} \otimes \lambda^{2,0} \rightr \\
& = & \leftr \lambda^{1,0} \otimes \lambda^{2,0} \rightr \oplus
\leftr \lambda^{2,1} \rightr
\eeq
Il existe un sous-espace $U^1$, $U(n)$-invariant tel que
\[
\lambda^{1,0} \otimes \lambda^{2,0} \simeq \lambda^{3,0} \oplus U^1
\]
Alors
\begin{equation}
\Lambda^1 \otimes \leftr \lambda^{2,0} \rightr \simeq \leftr
\lambda^{3,0} \rightr \oplus \leftr U^1 \rightr \oplus \leftr
\lambda^{2,1}_0 \rightr \oplus \Lambda^1 \label{W1234}
\end{equation}
est la d\'ecomposition en composantes irr\'eductibles de la repr\'esentation de $U(n)$ en dimension sup\'erieure ou \'egale \`a 6. Pour une d\'emonstration, on se reportera \`a \cite{gr4} en identifiant, avec les notations des auteurs
\[ W_1 \simeq \leftr
\lambda^{3,0} \rightr, \quad W_2 \simeq \leftr U^1 \rightr, \quad W_3
\simeq \leftr
\lambda^{2,1}_0 \rightr \quad \text{et} \quad W_4 \simeq \Lambda^1 \]
En dimension 4, $W_1$, $W_3$ sont r\'eduits \`a $\{0\}$.
La m\^eme d\'ecomposition est calcul\'ee par Falcitelli, Farinola, Salamon
\cite{fa} en se servant de l'algorithme de d\'ecomposition des produits
tensoriels expos\'e au chapitre 6 de \cite{sa} (figure 6.5). C'est
\'egalement la proc\'edure suivie dans cet article pour obtenir les
d\'ecompositions donn\'ees sans d\'emonstration \`a la section 4 (la preuve compl\`ete est d\'evelopp\'ee au chapitre 2, section 3 de \cite{bu3}).

Au lieu de demander que la premi\`ere d\'eriv\'ee de $\omega$ satisfasse
certaines conditions, on demandera, dans ce langage, que $\n \omega$
prenne ses valeurs dans certains sous-espaces invariants d\'efinis par
(\ref{W1234}) :
\begin{defi}
Soit $I \subset \{1,2,3,4\}$. On appelle vari\'et\'e de type $\bigoplus_{i
  \in I} W_i$, une vari\'et\'e presque hermitienne $(M,g,J)$ telle que la
d\'eriv\'ee covariante pour la connexion de Levi-Civita de la forme de
K\"ahler $\n \omega$ est une section du fibr\'e $\bigoplus_{i \in I} {\bf W_i} \subset \LLambda^1 \otimes
\leftr \llambda^{2,0} \rightr$. De plus, on appelle classe $\bigoplus_{i \in
  I} W_i$ l'ensemble de ces vari\'et\'es.
\label{Gray Hervella}
\end{defi}

Pour information et pour exemple, la classe $W_3 \oplus W_4$ est la
classe des vari\'et\'es hermitiennes, et $W_2$ est la classe des vari\'et\'es
symplectiques. Comme on voit, l'intersection de ces deux classes est
form\'ee de vari\'et\'es v\'erifiant $\n \omega =0$ ce qui correspond \`a la
d\'efinition alternative des vari\'et\'es k\"ahleriennes : la diff\'erentielle
de la forme de K\"ahler et le tenseur de Nijenhuis sont nuls en m\^eme temps. On verra section 5 une interpr\'etation utile des composantes de $\n \omega$ en fonction de $N$, $d\omega$ (voir figure 1).

On d\'efinit la connexion
hermitienne canonique $\nb$ par
\begin{equation}
\nb_X = \n_X - \frac{1}{2}J(\n_X J)
\end{equation}
On v\'erifie que $\nb J=0$. En outre quel que soit $X \in TM$, 
\begin{equation}
\overline{\delta}_X = \n_X -
\nb_X = \frac{1}{2}J(\n_X J)
\label{n-nb}
\end{equation}
anticommute \`a $J$ \`a cause de (\ref{nJ}).
La connexion $\nb$ est en fait l'unique connexion hermitienne ayant
cette propri\'et\'e. En effet si $\nt$ est une autre connexion hermitienne,
$\nb_X - \nt_X$ commute \`a $J$ quel que soit $X \in TM$. On peut donc toujours
d\'ecomposer
\begin{equation}
\begin{array}{ccccc}
\n - \nt & = & (\n - \nb) & + & (\nb - \nt) \\
\LLambda^1 \otimes \mathfrak{so}(M) & = & \LLambda^1 \otimes
\mathfrak{u}(M)^{\perp} & \oplus & \LLambda^1 \otimes
\mathfrak{u}(M)
\end{array}
\label{n-nt}
\end{equation}
Cette propri\'et\'e caract\'erise la \emph{connexion intrins\`eque} de la
structure $U(n)$, d\'efinie pour toute $G$-structure, $G
\subset SO(m)$. Le tenseur $\overline \delta = \n -
\nb$ est appel\'e (par abus de langage) \emph{torsion intrins\`eque} de la
structure $U(n)$. En effet, pour une connexion m\'etrique $\nt$, l'application qui \`a $\widetilde \delta = \n - \nt$ associe la torsion $T$ :
\[ T(X,Y)= \widetilde \delta_X Y - \widetilde \delta_Y X \]
est un isomorphisme. Maintenant, $\overline \delta$ est envoy\'e sur $\n
\omega$ par un isomorphisme $U(n)$-invariant. On peut donc aussi
regarder les composantes de ce tenseur dans la d\'ecomposition
(\ref{W1234}). Par l\`a, la d\'ecomposition de Gray-Hervella peut-\^etre
g\'en\'eralis\'ee \`a toute $G$-structure. Les articles \cite{ch,ma}
seront cit\'es dans la suite pour le cas $G=SU(n)$ (le premier
seulement pour $G=SU(3)$). On mentionne \'egalement sur ce sujet l'article ant\'erieur \cite{bo} de Bor, Hern\'andez Lamoneda.

\part{}

\section{L'espace de twisteurs r\'eduit et sa structure presque complexe}

Soit $Z(n)$ l'ensemble des endomorphismes unitaires de carr\'e
$-1$ de $\CM^n$. Un \'el\'ement $J$ de $Z(n)$ a l'inverse
$J^{-1}= \, ^t \!\bar J$ mais aussi, puisque $J^2=-1$, $J^{-1}=-J$
d'o\`u $J + \, ^t \! \bar J = 0$. Par cette remarque, $Z(n)$ est
l'intersection $U(n) \cap \, \mathfrak{u}(n)$ du groupe unitaire
et de son alg\`ebre de Lie. La donn\'ee d'un tel endomorphisme
diagonalisable est \'equivalente \`a la donn\'ee d'un couple
de sous-espaces complexes orthogonaux $(F,F^{\perp})$, les
sous-espaces propres de $J$ pour les valeurs propres $i$ et $-i$. En
classant suivant la dimension $p$ de $F$, on voit que $Z(n)$ a
plusieurs composantes connexes, chacune isomorphe \`a une
Grassmanienne complexe $\mathcal G_p(\CM^n)$. Le groupe $U(n)$ agit
transitivement par
\[(g,J) \mapsto gJg^{-1} \]
sur chaque composante connexe qui s'identifie ainsi \`a l'espace
homog\`ene $U(n)/ U(p) \times U(n-p)$. En particulier la composante
correspondant \`a $p=1$ est isomorphe \`a $\CM P^{n-1}$. En
revanche, on \'elimine les composantes singuli\`eres correspondant \`a la
multiplication par $i$ et $-i$ de $\CM^n$. De plus $Z(n)$ est muni
d'une structure presque complexe canonique $U(n)$-invariante
not\'ee ${\J}_n$, donn\'ee par la multiplication \`a gauche
par $J$ sur chaque espace tangent $T_J Z(n)$, identifi\'e \`a l'ensemble des endomorphismes de $\mathfrak u(n)$ qui anticommutent \`a $J$ (cf \cite{ob}).

\begin{defi}
Soit $(M,g,J_0)$ une vari\'et\'e presque hermitienne de
dimension $m=2n$ (NB : par \'economie de notation, on appelle $J_0$, dans cette partie, la structure presque complexe de la vari\'et\'e). L'espace de
twisteurs r\'eduit $Z$ de $M$ est le fibr\'e associ\'e du
fibr\'e principal $U(M)$ pour l'action de $U(n)$
sur $Z(n)$.
\end{defi}

Cette d\'efinition prend place dans un cadre tr\`es g\'en\'eral. Il s'agit d'un cas particulier d'une construction naturelle d'espaces de twisteurs sur des $G$-vari\'et\'es utilis\'ee dans \cite{ber,ob} (voir aussi \cite{bur} pour le cas des espaces sym\'etriques).

Le sous-fibr\'e vertical $T^V Z$ (tangent aux
fibres) est muni d'un endomorphisme $\J^V$ de carr\'e $-1$
copi\'e sur $\J_n$.

On adopte les notations suivantes : $j$ d\'esigne toujours un
point de $Z$, $\pi : Z \to M$ est la projection canonique et
$x=\pi(j)$. Bien s{\^u}r on peut voir $j$ comme un endomorphisme
de carr\'e $-1$ de $T_xM$, commutant avec $(J_0)_x$ et
compatible avec $g_x$. Alors $\J^V_j$ est la multiplication \`a
gauche par $j$ dans $T_j^V Z$. D'un autre c\^ot\'e on notera par une
majuscule, $J$, une section de Z, globale ou locale,
c'est-\`a-dire une structure presque complexe sur $M$ ou un
ouvert de $M$.

Toute connexion sur $U(M)$, c'est-\`a-dire toute
connexion hermitienne $\nt$ d\'efinit une connexion sur $Z$. On
note $\widetilde H \subset TZ$ la distribution horizontale correspondante. Elle permet de
compl\'eter $\J^V$ en une structure presque complexe $\J$ sur $Z$ en demandant que la restriction de $\J$
\`a $\widetilde H_j$, pour tout $j \in Z$, soit le relev\'e de $j$ lui-m\^eme, vu
comme structure presque complexe sur $T_x M$.

Maintenant si on utilise une autre connexion $\widehat \n$, on
obtient une structure presque complexe a priori diff\'erente
$\J'$. On note
\[ \eta_X = \nt_X - \widehat \n_X, \]
pour tout $X \in TM$.

\begin{prop}
Deux connexions hermitiennes d\'efinissent la m\^eme structure presque
complexe sur l'espace de twisteurs r\'eduit $Z$ si et seulement si
leur diff\'erence v\'erifie
\begin{equation}
[\eta_{JX},J] = J[\eta_X,J] \label{eta}
\end{equation}
pour tout vecteur $X$ et toute section $J$ de $Z$, c'est \`a
dire toute structure presque complexe sur $M$ commutant avec $J_0$.
\label{2 connexions 1 J}
\end{prop}

\begin{proof}
Soit $U$ un vecteur tangent \`a $Z$ en $j$. On appelle $X =
\pi_{*}(U)$ sa projection sur $T_xM$. Alors il existe une section
locale $J$ de $Z$ telle que $J_x=j$, $J_{*}(X)=U$. Par d\'efinition
de la d\'eriv\'ee covariante
\[ U = J_{*}(X) = \nt_X J + \widetilde X \]
o\`u $\widetilde X$ d\'esigne le relev\'e horizontal de $X$ dans
$\widetilde H_j$. Cette d\'ecomposition est en outre la
d\'ecomposition de $U$ en sa partie verticale et sa partie
horizontale, qui sert \`a calculer $\J U$. Mais de
m\^eme $U = \widehat \n_X J + \widehat X$ o\`u $\widehat X$ est le relev\'e de $X$ dans $\widehat H$, l'espace horizontal associ\'e \`a $\widehat \n$. D'o\`u on d\'eduit premi\`erement que
\[ \widehat X = \widetilde X + \nt_X J - \widehat \n_X J = \widetilde X + [\eta_X,J] \]
Puis
\[ \J U = \J^V(\nt_X J) + \widetilde{JX},\qquad \J'U = \J^V(\widehat \n_X J) + \widehat{JX} \]
On trouve que deux connexions $\nt$ et $\widehat \n$ d\'efinissent la
m\^eme structure presque complexe si et seulement si leur
diff\'erence v\'erifie
\[ [\eta_{JX},J] = \J^V [\eta_X,J] \]
Ce n'est rien d'autre que (\ref{eta}), par d\'efinition de $\J^V$.
\end{proof}

On peut r\'e\'ecrire (\ref{eta}) sous une autre forme. Pour cela, on
complexifie l'espace tangent et tous ses espaces de tenseurs et on
\'etend $\eta$ \`a $T^{\CM} M$ par $\CM$-lin\'earit\'e. Comme on l'a vu dans les pr\'eliminaires, chaque vecteur $X \in T^{\CM}M$ peut
\^etre d\'ecompos\'e suivant les sous-espaces propres
$T^{1,0}$, $T^{0,1}$ de $J_0$ en chaque point : $X = X^{1,0} + X^{0,1}$. Il peut aussi \^etre d\'ecompos\'e suivant les sous-espaces propres de $J$. Soit $X^{1,0}_J=\frac{1}{2}(X-iJX)$, dans cette partie, la projection de $X$ sur $T_J^{1,0}$, le sous-fibr\'e des vecteurs de type (1,0) \emph{par rapport \`a $J$}. De m\^eme, on note $X^{0,1}_J=\frac{1}{2}(X+iJX)$ la projection sur $T_J^{0,1}$. La condition donn\'ee \`a la proposition \ref{2 connexions 1 J} peut alors \^etre remplac\'ee par
\begin{equation} (\eta_{X^{1,0}_J} Y^{1,0}_J)^{0,1}_J = 0 \label{eta'} \end{equation}
Cette derni\`ere \'equation est elle-m\^eme \'equivalente \`a
\[ (\eta_{X^{0,1}_J} Y^{0,1}_J)^{1,0}_J = 0 \]
car $\eta$ est un tenseur r\'eel.

\vs

Une question naturelle maintenant est l'int\'egrabilit\'e de
$\J$, en fixant une connexion, ou une classe d'\'equivalence de
connexions d\'efinie par (\ref{eta'}). B\'erard Bergery et Ochiai \cite{ber} ou O'Brian et Rawnsley \cite{ob} ont donn\'e la r\'eponse sous la forme d'\'equations portant sur la torsion $T$ et la courbure $\widetilde R$ de $\nt$.

\begin{prop}[B\'erard Bergery, Ochiai]
La structure presque complexe $\J$ de l'espace de twisteurs r\'eduit $Z$ associ\'ee \`a la
connexion hermitienne $\nt$ est int\'egrable si et seulement si
\begin{equation}
T(X^{1,0}_J,Y^{1,0}_J)^{0,1}_J = 0 \label{T}
\end{equation}
et
\begin{equation}
(\Rt_{X^{1,0}_J,Y^{1,0}_J}Z^{1,0}_J)^{0,1}_J = 0,
\label{R}
\end{equation}
pour tous vecteurs $X,Y,Z \in T^{\CM} M$ et toute structure presque complexe $J$
commutant avec $J_0$.
\label{O'Brian & Rawnsley}
\end{prop}

\begin{NB}
A noter la ressemblance entre (\ref{eta'}) ou (\ref{T}) et (\ref{N}). En effet (\ref{T}) correspond \`a l'annulation de la partie horizontale du tenseur de Nijenhuis $\mathcal N$ de $\J$. Quant \`a (\ref{R}) elle est li\'ee \`a sa partie verticale. En effet, si $U$, $V$ sont deux champs de vecteurs tangents \`a $Z$, la projection sur $T^V Z$ du crochet $[U,V]$ au point $j$ ne d\'epend que de $U_j$, $V_j$ et vaut $[\widetilde R_{X,Y},j]$, o\`u $X = \pi_*(U)$ et $Y=\pi_*(V)$. 
\end{NB}

\section{R\'esolution des \'equations portant sur la torsion}

On s'attache premi\`erement \`a r\'esoudre l'\'equation
portant sur la torsion, c'est-\`a-dire qu'on cherche des
conditions {\it n\'ec\'essaires} pour que, une connexion
hermitienne $\nt$ \'etant fix\'ee, la structure presque complexe
associ\'ee $\J$ soit int\'egrable.

On a vu que $J$ agit en tant que section du fibr\'e adjoint $\mathfrak{so}(m)$ sur les espaces de tenseurs, en particulier sur l'espace des tenseurs de type torsion par
\[ J.T(X,Y) = JT(X,Y) -T(JX,Y) -T(X,JY) \]
Comme $J$ est diagonalisable, avec les valeurs propres $i$ et $-i$ sur $TM$, il l'est encore sur $\LLambda^2 \otimes TM$ avec les valeurs propres $3i$ (c'est \`a dire $i+i+i$), $-3i$,
$i$($=i+i-i$) et $-i$. Soit $T$ un vecteur propre pour la valeur
propre $-3i$.
\[ JT(X,Y) -T(JX,Y) -T(X,JY) = -3iT (X,Y) \]
On voit que $JT(X^{1,0}_J,Y^{1,0}_J)=-iT(X^{1,0}_J,Y^{1,0}_J)$ c'est-\`a-dire $T(X^{1,0}_J,Y^{1,0}_J) \in T^{0,1}_J$. De m\^eme on peut montrer que $T(X^{1,0}_J,Y^{1,0}_J) \in T^{1,0}_J$, si $T$ est un vecteur propre pour la valeur propre $-i$, et $T(X^{1,0}_J,Y^{1,0}_J)=0$ si $T$ est un vecteur propre pour les valeurs propres $3i$ ou $i$. Par cons\'equent, pour un tenseur quelconque $T$, (\ref{T}) est \'equivalente \`a l'annulation de la composante de $T$ suivant le sous-espace propre de $J$ associ\'e \`a $-3i$. De m\^eme, l'\'equation conjugu\'ee,
\begin{equation}
T(X^{0,1}_J,Y^{0,1}_J)^{1,0}_J = 0 \label{T'}
\end{equation}
correspond \`a l'annulation de la composante dans le sous-espace propre pour la valeur propre $3i$. Finalement (\ref{T}) et (\ref{T'}) sont les \'equations de la somme
directe des sous-espaces propres pour les valeurs propres $i$ et $-i$.

Maintenant, pour que la partie horizontale de $\mathcal N$ soit
nulle, par la proposition \ref{O'Brian & Rawnsley}, $T$ doit v\'erifier ces \'equations quel que soit $J$
commutant avec $J_0$ (pour un tenseur r\'eel en effet, l'une de ces \'equations implique imm\'ediatement l'autre). Appelons $\mathcal T$ l'espace des tenseurs de type
torsion qui le v\'erifient. Il est $U(n)$-invariant.
En effet, $J$ n'est pas en g\'en\'eral invariant par un
\'el\'ement $g \in U(n)$ mais envoy\'e par lui sur une autre
structure presque complexe $gJg^{-1}$ commutant avec $J_0$. C'est
pourquoi, pour calculer $\mathcal T$, on d\'ecompose
l'espace des tenseurs de type torsion en chaque point en
composantes irr\'eductibles de la repr\'esentation de $U(n)$
puis on fixe une structure presque complexe $J$ commutant avec
$J_0$. Si une composante $W$ a une intersection non nulle avec les
sous-espaces propres de $J$ correspondant aux valeurs propres $3i$ ou
$-3i$, elle n'appara{\^\i}t pas dans la d\'ecomposition de
$\mathcal T$. Dans le cas contraire, si $J$ n'a que les valeurs
propres $i$ et $-i$ sur $W$, par invariance, c'est aussi le cas de
toutes les autres structures presque complexes obtenues \`a
partir de $J$ par conjugaison avec un certain $g \in U(n)$,
c'est-\`a-dire appartenant \`a la m\^eme composante connexe
de l'espace de twisteurs $Z$. Pour obtenir une condition
d'int\'egrabilit\'e de $\J$, il suffit par cons\'equent de consid\'erer une section $J$ par composante connexe. Autrement, on peut chercher des conditions
d'int\'egrabilit\'e \emph{partielle}, en se restreignant
\`a une sous-vari\'et\'e connexe de $Z$. Les conditions obtenues sont a priori diff\'erentes pour des composantes connexes distinctes, comme on va montrer que c'est le
cas en dimension 8.

Comme on a vu dans les pr\'eliminaires, il revient au m\^eme de travailler
avec $T$ ou $\widetilde \delta = \n - \nt$. On choisit le second. C'est une section de $\LLambda^1 \otimes \mathfrak{so}(M) \simeq
\LLambda^1 \otimes \LLambda^2$.
D\'ecomposons cet espace sous l'action de $U(n)$.
\[ \Lambda^1 \otimes \Lambda^{2} = (\Lambda^1 \otimes \leftr
\lambda^{2,0} \rightr) \oplus (\Lambda^1 \otimes [\lambda^{1,1}] )
\]
Le premier sous-espace se d\'ecompose comme en (\ref{W1234}). Quant au second, il existe un sous-espace $U(n)$-invariant $U^2$ tel que
\[
\lambda^{1,0} \otimes \lambda^{1,1} \simeq \lambda^{2,1} \oplus  U^2
\]
Alors, en notant $U^2_0$ la composante primitive de $U^2$, on obtient la d\'ecomposition en composantes
irr\'eductibles :
\[
\Lambda^1 \otimes [ \lambda^{1,1} ] \simeq \Lambda^1 \oplus \leftr
\lambda^{2,1}_0 \rightr \oplus \Lambda^1 \oplus  \leftr U^2_0 \rightr
\]

\begin{prop} \label{CN}
Si $M$ est de dimension sup\'erieure \`a 10,
$\mathcal T \simeq 3\LLambda^1$.
Si $M$ est de dimension 6, $\mathcal T \simeq 3\LLambda^1 \oplus \leftr \llambda^{3,0} \rightr$.
\label{espace torsions}
\end{prop}

\begin{proof}
Soit $J$ un endomorphisme de carr\'e $-1$ commutant avec $J_0$. On
diagonalise simultan\'ement $J$ et $J_0$ en chaque point. Soient
$F$ le sous-fibr\'e de $T^{\CM} M$ o\`u $J$ co{\"\i}ncide avec
$J_0$, $G$ le sous-fibr\'e sur lequel il vaut $-J_0$. On a $F \oplus
G = T^{\CM}M$ et m\^eme
\[ \llambda^{1,0} =  \llambda^{1,0}F \oplus  \llambda^{1,0}G, \quad \llambda^{0,1} =  \llambda^{0,1}F \oplus  \llambda^{0,1}G
\]
o\`u $\llambda^{1,0}F$ est l'ensemble des formes de type (1,0) qui s'annulent sur $G$, etc. Une autre fa\c con de
le dire est $\llambda^{1,0}F =
\llambda^{1,0} \cap {\bf l}^{1,0}$, o\`u ${\bf l}^{1,0}$ est l'espace des
1-formes de type $(1,0)$ {\it par rapport \`a $J$}, $\llambda^{1,0}G
= \llambda^{1,0} \cap {\bf l}^{0,1}$, etc. On a donc aussi
\beq
{\bf l}^{1,0} & = & {\bf l}^{1,0}F \oplus {\bf l}^{1,0}G \\
        & = & \llambda^{1,0}F \oplus \llambda^{0,1}G
\eeq
o\`u ${\bf l}^{1,0}F$ est comme on s'y attend l'ensemble des formes de type (1,0) {\it par rapport \`a $J$} qui s'annulent sur $G$, etc.

Par la discussion pr\'ec\'edente, les composantes de $\mathcal T$ sont celles qui ne rencontrent pas
$\bo{3} {\bf l}^{1,0}$ ni $\bo{3} {\bf l}^{0,1}$.
\[ \bo{3} {\bf l}^{1,0} = \bo{3} \llambda^{1,0}F \oplus
(\bo{2} \llambda^{1,0}F \otimes \llambda^{0,1}G) \oplus
(\llambda^{1,0}F \otimes \bo{2} \llambda^{0,1}G) \oplus \bo{3}
\llambda^{0,1}G
\]
Examinons premi\`erement l'intersection de chaque terme avec $\llambda^1 \otimes \llambda^{2,0}$.

Le premier, $\bo{3} \llambda^{1,0}F$, a une intersection non r\'eduite
\`a $\{0\}$ avec $\llambda^{3,0}$: $\llambda^{3,0}F$. Il a aussi une intersection non
nulle avec ${\bf U}^1$ qu'on peut noter ${\bf U}^1 F$.

Le deuxi\`eme, $\bo{2} \llambda^{1,0}F \otimes \llambda^{0,1}G$, a une
intersection non nulle avec $\llambda^{2,1}_0$ : $\llambda^{2,0}F \otimes
\llambda^{0,1}G$. (Il est \`a noter que, comme $F$ et $G$ sont
$J$-stables, orthogonaux, $\llambda^{2,0}F \otimes \llambda^{0,1}G$ a une
intersection nulle avec $\llambda^{1,0}$. En effet,
celui-ci est obtenu dans $\llambda^{2,1}$ en faisant le produit
ext\'erieur avec $\omega \in \llambda^{1,1}F \oplus \llambda^{1,1}G$.)

Enfin les deux derniers sous-espaces $\llambda^{1,0}F \otimes
\bo{2} \llambda^{0,1}G$ et $\bo{3} \llambda^{0,1}G$ ont une
intersection nulle avec $\llambda^1 \otimes \llambda^{2,0}$.

On proc\`ede de m\^eme pour les composantes de $\llambda^1
\otimes \llambda^{1,1}$ : le sous-fibr\'e isomorphe \`a
$\llambda^{2,1}_0$ a une intersection non nulle avec $\bo{3}
{\bf l}^{1,0}$: $\llambda^{2,0}F \otimes \llambda^{0,1}G$, quant \`a ${\bf U}^2_0$, il
intersecte $\bo{3} {\bf l}^{1,0}$ suivant $\boldsymbol{\sigma}^{2,0}F \otimes
\llambda^{0,1}G$. (NB : On note $\sigma^{p,q}=\bigodot^p
\lambda^{1,0} \otimes \bigodot^q \lambda^{0,1}$. On voit par
exemple que $U^2 = \sigma^{2,1}$.)

Finalement il semble que tous les sous-espaces irr\'eductibles de
$\LLambda^1 \otimes \LLambda^2$, sauf ceux isomorphes \`a $\LLambda^1$,
rencontrent le sous-espace propre de $J$ pour la valeur propre
$3i$. Cependant il faut voir que si les
dimensions de $F$ ou $G$ sont petites, certaines des intersections pr\'ec\'edentes sont
r\'eduites \`a $\{0\}$. Ainsi, si $F$ est de dimension $2$ ou $4$,
$\llambda^{3,0}F = \{0\}$. Mais alors, pourvu que $M$ soit de
dimension sup\'erieure \`a 10, $G$ est de dimension sup\'erieure \`a
$6$ et $\llambda^{3,0}$ rencontre $\bo{3} {\bf l}^{0,1}$, s'il ne rencontre
pas $\bo{3} {\bf l}^{1,0}$.
\end{proof}

\begin{coro}
Soit $(M,g,J_0)$ une vari\'et\'e presque hermitienne de
dimension $m=2n$. Une condition n\'ec\'essaire d'int\'egrabilit\'e de $\J$ est que $M$ soit localement conform\'ement k{\"a}hlerienne en dimension sup\'erieure \`a 10 ou de
type $W_1 \oplus W_4$ en dimension 6.
\label{W14}
\end{coro}

\begin{proof}
Comme on le voit dans la d\'ecomposition (\ref{n-nt}), la partie de $\widetilde \delta$ dans $\LLambda^1 \otimes
\mathfrak{u}(M)^{\perp} \simeq \LLambda^1 \otimes \leftr \llambda^{2,0}
\rightr$ est fixe, \'egale
\`a $\overline \delta = \n - \nb$. D'apr\`es une remarque pr\'ec\'edente,
celle-ci (\ref{n-nb}), au lieu de $\n \omega$, peut servir \`a d\'efinir
le type de la vari\'et\'e presque hermitienne. Les conditions impos\'ees \`a
$T$ ou $\widetilde \delta$, proposition \ref{espace torsions}, impliquent
que les composantes suivant ${\bf W}_2 \simeq \leftr {\bf U}^1
\rightr$, ${\bf W}_3 \simeq \leftr \llambda^{2,1}_0 \rightr$ -- ainsi que ${\bf W}_1 \simeq \leftr \lambda^{3,0} \rightr$, en dimension 6 -- sont nulles.
\end{proof}

La particularit\'e de la dimension 6 vient de ce que $dim_{\CM}F$ et
$dim_{\CM}G$ sont toutes deux strictement inf\'erieures \`a $6$ pourvu que
$J$ ne soit pas \'egal \`a $J_0$, ni \`a son oppos\'e. Alors, $\llambda^{1,0}F$,
$\llambda^{0,1}G$ sont de dimension 1 ou 2 et $\llambda^{3,0}F$,
$\llambda^{0,3}G$ sont r\'eduits \`a $\{0\}$. En dimension 8,
c'est encore le cas si $dim_{\CM}F = dim_{\CM}G = 4$, c'est \`a dire si
$J$, $J_0$ d\'eterminent la m\^eme orientation de $M$. Appelons $Z_1$ le
sous-fibr\'e de $Z$ engendr\'e par ces sections,
$Z_2$ le sous-fibr\'e suppl\'ementaire. Comme on a vu, les composantes
connexes de $Z$ sont caract\'eris\'ees par les dimensions des sous-espaces $F$, $G$.
On en d\'eduit que $Z_1$ est connexe tandis que $Z_2$ a deux
composantes connexes correspondant respectivement \`a $dim_{\CM} F=2$,
$dim_{\CM} G=6$ et $dim_{\CM} F=6$, $dim_{\CM} G=2$. Si
$J$ est une section de $Z_2$, un parmi $\llambda^{3,0}F$,
$\llambda^{3,0}G$, selon les cas, n'est pas r\'eduit \`a $\{0\}$. Par cons\'equent,
quel que soit le
tenseur $T$ de $\leftr \llambda^{3,0} \rightr$, il existe $J \in Z_2$
tel que $T$ ne v\'erifie pas (\ref{T}). Ce qui nous int\'eresse alors, en
cette dimension, est moins le tenseur $\J$ de l'espace de twisteurs r\'eduit que ses restrictions \`a $Z_1$, $Z_2$, not\'ees $\J_1$, $\J_2$. Pour une vari\'et\'e $M$ de dimension 8, on note $\mathcal T_1$ (resp. $\mathcal
T_2$) le sous-espace de l'espace des tenseurs de torsion
abstraits dont les \'el\'ements v\'erifient (\ref{T}), quelle que soit la
section $J$ de $Z_1$ (resp. $Z_2$). Autrement dit, $\mathcal T_1$ (resp. $\mathcal T_2$) d\'etermine une
condition n\'ec\'essaire d'int\'egrabilit\'e de $\J_1$ (resp. $\J_2$).

\begin{prop}
En dimension 8, $\mathcal T_1 \simeq 3\LLambda^1 \oplus \leftr \llambda^{3,0} \rightr$ mais
$\mathcal T_2 \simeq 3\LLambda^1$.
\end{prop}

\begin{coro}
Soit $M$ une vari\'et\'e presque hermitienne de dimension 8. La structure presque complexe
$\J_1$ de $Z_1$ est int\'egrable seulement si la
vari\'et\'e est de type $W_1 \oplus W_4$ mais pour que $\J_2$ soit
int\'egrable il faut que $M$ soit localement conform\'ement
k{\"a}hlerienne.
\end{coro}

\vs

Aucun travail suppl\'ementaire n'est n\'ec\'essaire pour
r\'eduire ($\ref{eta}$), qui est identique \`a ($\ref{T}$)
\`a ceci pr\`es que $\eta$ vit dans le sous-espace $\LLambda^1
\otimes [ \llambda^{1,1}]$ de $\LLambda^1 \otimes \LLambda^2$. On en
d\'eduit que deux connexions hermitiennes d\'efinissent la
m\^eme structure presque complexe de l'espace des twisteurs
r\'eduits si et seulement si leur diff\'erence est dans
$\LLambda^1 \oplus \LLambda^1$. On d\'efinit de cette fa{\c c}on
une relation d'\'equivalence sur les connexions hermitiennes de
$(M,g,J_0)$ de telle sorte que la construction de $\J$ sur
l'espace de twisteurs r\'eduits est une injection de l'ensemble
des classes d'\'equivalence dans l'espace des structures presque
complexes de $Z$. L'\'enonc\'e suivant regroupe les deux
r\'esultats :

\begin{prop}
La structure presque complexe $\J$ de l'espace de twisteurs r\'eduit associ\'ee \`a une
connexion hermitienne $\nt$ ne peut \^etre int\'egrable que si cette derni\`ere appartient \`a la m\^eme classe d'\'equivalence que $\nb$, la connexion hermitienne canonique, i.e. $\nt$, $\nb$ d\'efinissent la m\^eme structure presque complexe sur $Z$. De plus la torsion
de $\nb$ doit \^etre une section de $\LLambda^1 \simeq {\bf W}_4$ en dimension sup\'erieure \`a 8 ou de
$\LLambda^1 \oplus \leftr \llambda^{3,0} \rightr \simeq {\bf W}_1 \oplus
{\bf W}_4$ en dimension 6.
\end{prop}

On est donc fond\'e dans la suite (recherche de conditions
suffisantes) \`a s'int\'eresser seulement \`a la connexion
hermitienne canonique $\nb$.

\section{Conditions d'int\'egrabilit\'e}

On veut maintenant r\'esoudre l'\'equation (\ref{R})
portant sur la courbure. Soit $\Rb$ la courbure de $\nb$.
De m\^eme que pour la torsion,
$\mathcal T$ \'etait le plus gros sous-espace $U(n)$-invariant
tel qu'une structure presque complexe $J$ commutant avec $J_0$,
agissant sur les tenseurs de type torsion, n'y ait que les valeurs
propres $i$ et $-i$ (jamais $3i$ ni $-3i$), de m\^eme $\mathcal
R$, l'ensemble des tenseurs qui v\'erifient ($\ref{R}$), est
l'espace des tenseurs de courbure hermitienne
abstraits, isomorphe \`a $\LLambda^2 \otimes
[ \llambda^{1,1}]$, \`a l'exclusion des composantes irr\'eductibles qui
admettent $\pm 4i$, comme valeur propre de $J$.

On d\'ecompose $\Lambda^2 \otimes [ \lambda^{1,1} ]$ sous l'action de $U(n)$.
\[ \lambda^2 \otimes \lambda^{1,1} = (\lambda^{2,0} \otimes
\lambda^{1,1}) \oplus (\lambda^{0,2} \otimes \lambda^{1,1}) \oplus
(\lambda^{1,1} \otimes \lambda^{1,1})
\]
Il existe un sous-espace $V^1$ tel que
\[ \lambda^{1,1} \otimes \lambda^{1,1} = \lambda^{2,2} \oplus V^1
\oplus \overline{V^1} \oplus \sigma^{2,2}
\]
Chaque terme de cette somme se d\'ecompose encore en :
\[ \lambda^{2,2} \simeq \lambda^{2,2}_0 \oplus \lambda^{1,1}_0 \oplus \RM
\]
\[ V^1 \simeq V^1_0 \oplus \lambda^{1,1}_0
\]
\[ \sigma^{2,2} \simeq \sigma^{2,2}_0 \oplus \lambda^{1,1}_0 \oplus \RM
\]
On rappelle que $\lambda^{p,q}_0$ est d\'efinie comme la repr\'esentation de $U(n)$ de
poids dominant, dans les coordonn\'ees standard,
\[
(\underbrace{1,\ldots,1}_p,0,\ldots,0,\underbrace{-1,\ldots,-1}_q)
\]
D'autre part, le poids dominant de $\sigma^{p,q}_0$ est $(p,0,\ldots,0,-q)$. En
particulier $\lambda^{1,1}_0=\sigma^{1,1}_0$ est de poids dominant
$(1,0,...,0,-1)$ et $\sigma^{2,2}_0$, de poids dominant
$(2,0,\ldots,0,-2)$, est le produit de Cartan de $\lambda^{1,1}_0$
avec elle-m\^eme. Maintenant, $U^1$ d\'efini pr\'ec\'edemment est la
repr\'esentation irr\'eductible de poids dominant
$(2,0,\ldots,-1)$. Enfin $V^1_0$ est d\'efinie comme la
repr\'esentation irr\'eductible de $U(n)$ de poids dominant
$(2,0,\dots,0,-1,-1)$. A noter que dans le produit sym\'etrique
il ne reste que
\[ \lambda^{1,1} \odot \lambda^{1,1} = \lambda^{2,2} \oplus \sigma^{2,2}
\]
et dans la dimension qui nous int\'eresse tout particuli\`erement, la
dimension 6 :
\[ [\lambda^{1,1}] \otimes [\lambda^{1,1}] = 3[\lambda^{1,1}_0] \oplus
2\RM \oplus
\leftr V^1_0 \rightr \oplus [\sigma^{2,2}_0]
\]

Il existe aussi $V^2$ tel que
\beq
\lambda^{2,0} \otimes \lambda^{1,1} & = & \lambda^{3,1} \oplus V^2
\\
& = & \lambda^{3,1}_0 \oplus 2\lambda^{2,0} \oplus V^2_0 \oplus
\sigma^{2,0}
\eeq
o\`u $V^2_0$ est la repr\'esentation irr\'eductible de poids dominant $(2,1,0
\ldots,0,-1)$. En dimension 6, en revenant aux espaces de
repr\'esentation r\'eels :
\[
\leftr \lambda^{2,0} \rightr \otimes [\lambda^{1,1}] = \leftr
\lambda^{2,0} \rightr \oplus \leftr V^2_0 \rightr \oplus
\leftr \lambda^{2,0} \rightr \oplus \leftr \sigma^{2,0} \rightr
\]
Le sous-espace propre de $J$ pour la valeur $4i$ dans
$\bo{4}\LLambda^1_{\CM}$ est
\begin{multline*}
\bo{4}{\bf l}^{1,0} = \bo{4}\llambda^{1,0}F \oplus (\bo{3}\llambda^{1,0}F
\otimes \llambda^{0,1}G) \oplus (\bo{2}\llambda^{1,0}F \otimes
\bo{2}\llambda^{0,1}G) \\ \oplus (\llambda^{1,0}F \otimes
\bo{3}\llambda^{0,1}G) \oplus \bo{4}\llambda^{0,1}G
\end{multline*}

Or,

-- $\bo{4}\llambda^{1,0}F$ n'a d'intersection non nulle qu'avec
$\llambda^{2,0} \otimes \llambda^{2,0}$. De m\^eme
$\bo{4}\llambda^{0,1}G$ ne nous int\'eresse pas.

-- $\bo{3}\llambda^{1,0}F \otimes \llambda^{0,1}G$ a une
intersection non nulle avec $\llambda^{3,1}_0$ (\'egale \`a
$\llambda^{3,0}F \otimes \llambda^{0,1}G$) et avec ${\bf V}^2$ (isomorphe
\`a ${\bf U}^1 F \otimes \llambda^{0,1}G$).

-- $\llambda^{1,0}F \otimes \bo{3}\llambda^{0,1}G$ est le conjugu\'e
du pr\'ec\'edent.

-- Enfin $\bo{2}\llambda^{1,0}F \otimes \bo{2}\llambda^{0,1}G$ a une
intersection non nulle avec $\llambda^{2,2}_0$ ($\llambda^{2,0}F
\otimes \llambda^{0,2}G$), ${\bf V}^1_0$ ($\llambda^{2,0}F \otimes
\boldsymbol{\sigma}^{0,2}G$), $\overline{{\bf V}^1_0}$ et $\boldsymbol{\sigma}^{2,2}_0$
($\boldsymbol{\sigma}^{2,0}F \otimes \boldsymbol{\sigma}^{0,2}G$).

Par cons\'equent les composantes qu'il faut \'eliminer sont $[
\llambda^{2,2}_0 ]$, $\leftr {\bf V}^1_0 \rightr$, $[ \boldsymbol{\sigma}^{2,2}_0 ]$,
$\leftr {\bf V}^2_0 \rightr$, $\leftr \llambda^{3,1}_0 \rightr$, en
dimension sup\'erieure \`a 10 et $\leftr {\bf V}^1_0 \rightr$, $[
\boldsymbol{\sigma}^{2,2}_0 ]$, $\leftr {\bf V}^2_0 \rightr$, en dimension 6.

Cependant,
\begin{prop}
Soit $(M,g,J)$ une vari\'et\'e presque hermitienne de type $W_1 \oplus
W_4$, de dimension 6. Les composantes dans $\leftr {\bf V}^1_0 \rightr$ et
$\leftr {\bf V}^2_0 \rightr$ de la courbure $\Rb$ de la connexion hermitienne
canonique $\nb$ sont nulles.
\end{prop}

\begin{proof}
Pour le voir on d\'ecompose
\begin{equation}
\overline \delta = \xi + \vartheta
\label{xi+theta}
\end{equation}
o\`u $\xi$ et $\vartheta$ sont des sections de ${\bf W}_1$ et ${\bf W}_4 \subset
\LLambda^1 \otimes \leftr \llambda^{2,0} \rightr$, respectivement. D\`es
lors $\Rb$ est li\'ee \`a la courbure riemannienne par
\begin{eqnarray}
\nonumber \Rb_{X,Y} & = & R_{X,Y} - (\nb_X \xi)_Y + (\nb_Y \xi)_X \\
& & - (\nb_X \vartheta)_Y + (\nb_Y \vartheta)_X - \vartheta_{\vartheta_X Y -
  \vartheta Y X} - [\vartheta_X,\vartheta_Y] \label{Rb=R} \\
\nonumber & & - \xi_{\xi_X Y - \xi Y X}  - [\xi_X,\xi_Y] \\
\nonumber & & - \vartheta_{\xi_X Y -
  \xi Y X} - \xi_{\vartheta_X Y - \vartheta Y X} - [\xi_X,\vartheta_Y]
- [\xi_Y,\vartheta_X]
\end{eqnarray}
Comme $\nb$ est une connexion hermitienne, $\nb \xi$ vit en tout
point dans un espace isomorphe, comme espace de repr\'esentation
de $U(n)$, \`a $\Lambda^1 \otimes \leftr \lambda^{3,0} \rightr$
; c'est aussi le cas des tenseurs aparaissant \`a la
derni\`ere ligne car $W_4 \simeq \Lambda^1$ ; les tenseurs
aparaissant \`a la deuxi\`eme ligne vivent dans un espace
isomorphe \`a $\LLambda^1 \otimes \LLambda^1$ ; enfin ceux
aparaissant \`a la troisi\`eme ligne vivent dans $\leftr
\llambda^{3,0} \rightr \otimes \leftr \llambda^{3,0} \rightr$. Or
cette \'equation permet de d\'ecomposer la courbure
hermitienne (en utilisant le projecteur $b$ d\'efini par la
permutation circulaire sur les trois premi\`eres variables) en
\begin{equation}
\Rb = \Rb_0 + \Rb_1
\label{R0+R1}
\end{equation}
o\`u $\Rb_0$ v\'erifie l'identit\'e de Bianchi, i.e. est un
tenseur de type courbure k{\"a}hlerienne, et $\Rb_1$ ne d\'epend
que de $b(\Rb)$, autrement dit, puisque $b(R)=0$, uniquement
des autres tenseurs aparaissant dans (\ref{Rb=R}). L'espace des
tenseurs de courbure k{\"a}hlerienne est
$\mathcal{K}=[\boldsymbol{\sigma}^{2,2}_0] \oplus [\llambda^{1,1}_0] \oplus
{\bf R}$. D\'ecomposons les autres espaces de repr\'esentation de
$U(n)$ en composantes irr\'eductibles.
\begin{equation}
\Lambda^1 \otimes \Lambda^1 = [ \lambda^{1,1}_0 ] \oplus \leftr \lambda^{2,0} \rightr \oplus \leftr \sigma^{2,0}
\rightr \oplus 2\RM
\label{L1xL1}
\end{equation}
\[
\Lambda^1 \otimes \leftr \lambda^{3,0} \rightr = \leftr V^3 \rightr
\oplus \leftr \lambda^{4,0} \rightr \oplus \leftr \lambda^{3,1}_0
\rightr \oplus \leftr \lambda^{2,0} \rightr
\]
o\`u $V^3$ est la repr\'esentation irr\'eductible (complexe) de poids
dominant $(2,1,1,0,\ldots,0)$. En dimension 6, $\leftr \lambda^ {4,0}
\rightr$ et $\leftr \lambda^{3,1}_0 \rightr$ n'aparaissent plus dans
cette d\'ecomposition. Enfin, en dimension 6,
\[ \leftr \lambda^{3,0} \rightr \otimes
\leftr \lambda^{3,0} \rightr = \leftr \lambda^{3,0} \rightr \oplus
\RM \oplus \RM
\]
Par cons\'equent ces trois espaces, et $\mathcal K$, ont une
intersection r\'eduite \`a $\{0\}$ avec $\leftr {\bf V}^1_0 \rightr$ et $\leftr
{\bf V}^2_0 \rightr$.
\end{proof}

Par ailleurs seul $R$, dans (\ref{Rb=R}), a \'eventuellement une
composante non nulle dans $[\boldsymbol{\sigma}^{2,2}_0 ]$ donc la condition
porte indiff\'eremment sur la courbure hermitienne ou la
courbure riemannienne. C'est ce qui autorise l'interpr\'etation de
O'Brian et Rawnsley dans \cite{ob} en termes de tenseur de Bochner
g\'en\'eralis\'e tel que d\'efini dans \cite{tr} pour la courbure
riemannienne de toute vari\'et\'e presque hermitienne.

\begin{theo} \label{dim6}
Soit $M$ une vari\'et\'e presque hermitienne de dimension 6 et
$\nt$ une connexion
hermitienne. La structure presque complexe $\J$ de l'espace de
twisteurs r\'eduit associ\'ee \`a $\nt$ est int\'egrable si et
seulement si les trois conditions suivantes sont satisfaites \\
(\romannumeral 1) $\nt$ d\'efinit la m\^eme structure presque complexe
sur $Z$ que la connexion hermitienne canonique $\nb$. \\
(\romannumeral 2) La vari\'et\'e est de type $W_1 \oplus W_4$. \\
(\romannumeral 3) Le tenseur de Bochner (ou la composante de $R$ dans $[\boldsymbol{\sigma}^{2,2}_0]$) est nul.
\end{theo}

En particulier,
\begin{coro}
L'unique vari\'et\'e strictement NK de dimension 6, compacte, simplement connexe, telle que l'espace de
twisteurs r\'eduit est muni d'une structure presque complexe $\J$ int\'egrable est la sph\`ere $S^6$.
\label{S6}
\end{coro}

\begin{proof}
Une vari\'et\'e NK est une vari\'et\'e de type $W_1$ dans la classification de Gray Hervella. De fa\c con \'equivalente $\n \omega$ ou $\overline \delta$ sont totalement antisym\'etriques. Cela correspond \`a $\vartheta = 0$ dans (\ref{xi+theta}). En outre, la torsion de la connexion hermitienne canonique est parall\`ele :
\begin{equation}
\nb \xi = 0,
\label{nbxi}
\end{equation}
comme d\'emontr\'e en toute g\'en\'eralit\'e dans \cite{ki}. D\`es lors (\ref{Rb=R}) se simplifie grandement. Soit $(\xi^2)$ le tenseur de $\LLambda^2 \otimes \mathfrak{so}(M)$
d\'efini par $(\xi^2)_{X,Y} = \xi_{\xi_X Y - \xi_Y X} +
[\xi_X,\xi_Y]$ :
\[ \Rb = R - (\xi^2)
\]
De plus, en dimension 6, les vari\'et\'es NK non k\"ahleriennes admettent une r\'eduction naturelle \`a $SU(3)$ associ\'ee \`a $d\omega$ ou $\xi$ (voir section 5). Alors (\ref{nbxi}) implique que $\nb$ est
une connexion $SU(3)$ et on peut d\'ecomposer $\Rb$ comme en (\ref{R0+R1}) o\`u cette fois $\Rb_0 \in \mathcal{K}(\mathfrak{su}(3))$ a les
propri\'et\'es alg\'ebriques du tenseur de courbure d'une
vari\'et\'e Calabi-Yau et
$\Rb_1$ ne d\'epend que de $(\xi^2)$. Finalement
\[ \Rb \in [\boldsymbol{\sigma}^{2,2}_0 ] \oplus {\bf R} \]
Telle est la forme extr\^emement simple de la courbure hermitienne
mais aussi riemannienne d'une vari\'et\'e strictement NK de
dimension 6.
\[ R = \Rb + (\xi^2) \in [\boldsymbol{\sigma}^{2,2}_0 ] \oplus {\bf R} \]
On retrouve de cette fa\c con que les vari\'et\'es NK non k\"ahleriennes, comme les vari\'et\'es de Calabi Yau en dimension 6
sont d'Einstein (voir \cite{gr3}). Maintenant, pour que $\J$ soit int\'egrable il faut, par le
th\'eor\`eme \ref{dim6}, (\romannumeral 3), que la composante de la courbure dans
$[\boldsymbol{\sigma}^{2,2}_0 ]$ soit nulle, c'est-\`a-dire que $R$ soit le
tenseur de courbure de la sph\`ere.
\end{proof}

En dimension sup\'erieure \`a 10, (\romannumeral 2) est remplac\'ee, conform\'ement au corollaire \ref{W14}, par : $M$ est localement conform\'ement k{\"a}hlerienne c'est-\`a-dire, cette fois,
$\xi=0$ dans (\ref{xi+theta}). Par cons\'equent,
compte-tenu de (\ref{L1xL1}), la m\^eme remarque
pr\'ec\'edant le th\'eor\`eme \ref{dim6} est valable, bien
qu'on ne soit plus en dimension 6. Comme en outre le tenseur de
Bochner est un invariant conforme on a

\begin{theo}
Soit $M$ une vari\'et\'e presque hermitienne de dimension
sup\'erieure \`a 10. La structure presque complexe $\J$ de
l'espace de twisteurs r\'eduit, associ\'e \`a une connexion
hermitienne $\nt$ est
int\'egrable si et seulement si \\
(\romannumeral 1) $\nt$ d\'efinit la m\^eme structure presque
complexe que $\nb$. \\
(\romannumeral 2) $M$ est localement conforme \`a une vari\'et\'e
k{\"a}hlerienne, Bochner-plate.
\end{theo}

Il s'agit du r\'esultat obtenu par O'Brian et Rawnsley \cite{ob}. Il reste
valable en dimension 8 \`a condition de se restreindre \`a la
sous-vari\'et\'e $Z_2$ de $Z$ d\'efinie plus haut.

\part{}

\section{Vari\'et\'es de type $W_1 \oplus W_4$}

La question pos\'ee dans cette section est la suivante : toutes les les
vari\'et\'es de type $W_1 \oplus W_4$, de dimension 6,
sont-elles localement conformes \`a des vari\'et\'es NK ? Autrement dit, la forme de Lee $\theta$
repr\'esente-t-elle toujours un changement conforme local, c'est-\`a-dire est-elle ferm\'ee (localement exacte) ?

Commen\c cons par \'etablir quelques notations. Soit $(M,g,J)$ une vari\'et\'e pres\-que hermitienne de
dimension $m=2n$. On appelle $w_i$, $i=1,2,3,4$ les composantes de $\n \omega$ dans la d\'ecomposition (\ref{W1234}) servant \`a d\'efinir la classification de Gray, Hervella. L'op\'erateur de Bianchi (ou la permutation circulaire sur les trois variables) envoie $\LLambda^1 \otimes \leftr
\llambda^{2,0} \rightr$ surjectivement sur l'espace des 3-formes $\LLambda^3$ et $\n \omega$ sur $d \omega$. Il est invariant sous l'action de $U(n)$ et son noyau est ${\bf W}_2$ par cons\'equent :
\[ {\bf W}_1 \oplus {\bf W}_3
\oplus {\bf W}_4 \simeq \LLambda^3 \simeq \leftr \llambda^{3,0} \rightr \oplus \leftr
\llambda^{2,1}_0 \rightr \oplus \LLambda^1
\]
La derni\`ere composante est plus
pr\'ecis\'ement dans $\LLambda^3$ l'image de l'application qui
fait le produit ext\'erieur des 1-formes avec $\omega$. On propose d'appeler $\theta$ la 1-forme apparaissant dans la d\'ecomposition de $d\omega$. D'autre part, on note $\psi$ sa composante dans $\leftr \llambda^{3,0} \rightr$,  de sorte que $d\omega$ s'\'ecrit finalement :
\[
d\omega = \psi + (d\omega)^{2,1}_0 + \omega \wedge \theta
\]
En dimension 6, la donn\'ee
d'une 3-forme de type (3,0)+(0,3), de norme constante, est \'equivalente \`a la donn\'ee d'une 3-forme
volume complexe d\'eterminant une r\'eduction de la
vari\'et\'e presque hermitienne \`a $SU(3)$, comme on va le voir bient\^ot.

Bien que $\n \omega$, $\overline \delta$ soient des tenseurs \'equivalents, pour notre probl\`eme et du point de vue des repr\'esentations, on pr\'ef\`ere noter diff\'eremment, pour plus de clart\'e, leurs composantes dans la d\'ecomposition (\ref{W1234}). Ainsi, la composante de $\overline \delta$ dans le sous-fibr\'e de $\LLambda^1 \otimes \mathfrak{u}(M)^{\perp}$ isomorphe \`a $\leftr \llambda^{3,0}
\rightr$ et correspondant \`a $w_1$ ou $\psi$ sera not\'ee $\xi$ et $\vartheta$, comme \`a la section pr\'ec\'edente, la composante dans le sous-fibr\'e isomorphe \`a $\LLambda^1$, li\'ee \`a $w_4$ ou $\theta$ par
\begin{equation}
4\vartheta_X = X^{\flat} \wedge \theta - JX^{\flat} \wedge
J\theta \label{(theta)},
\end{equation}
o\`u $\flat$ d\'esigne l'\og isomorphisme musical \fg \ qui \`a
un vecteur associe sa 1-forme duale par la m\'etrique. (Notons
\'egalement que l'expression ci-dessus identifie un endomorphisme
antisym\'etrique de l'espace tangent $\vartheta_X$ et
une 2-forme diff\'erentielle.)

\begin{figure}[h]
\caption{Classes de Gray \& Hervella}
\begin{picture}(300,300)(-40,-40)
\put(0,0){\framebox(110,110){\shortstack{\large{$\mathbf{W_2
\simeq
      \leftr U^1 \rightr}$} \\ \\ $N_2$  }}}
\put(111,0){\framebox(110,110){\shortstack{\large{$\mathbf{W_4
\simeq
      \LLambda^1}$} \\ \\ $\vartheta$ ou $\theta$}}}
\put(0,110){\framebox(110,110){\shortstack{$\xi$ ou $\psi$ ou $N_1$ \\ \\
       \large{$\mathbf{W_1 \simeq \leftr \llambda^{3,0} \rightr}$}}}}
\put(111,110){\framebox(110,110){\shortstack{$(d\omega)^{1,2}_0$
      \\ \\ \large{$\mathbf{W_3 \simeq \leftr \llambda^{2,1}_0 \rightr}$}}}}
\put(100,165){\oval(220,90)[t]} \put(100,165){\oval(220,90)[bl]}
\put(210,80){\line(0,1){100}} \put(165,80){\oval(90,180)[b]}
\put(120,80){\line(0,1){40}} \put(100,120){\line(1,0){20}}
\put(8,-10){\dashbox{10}(94,240)} \put(8,30){\line(-3,-1){49}}
\put(-55,8){\makebox(10,10){\large{$N$}}}
\put(0,80){\line(-1,0){35}}
\put(-55,77){\makebox(10,10){\large{$\nabla \omega$}}}
\put(-10,170){\line(-3,1){29}}
\put(-55,177){\makebox(10,10){\large{$d\omega$}}}
\end{picture}
\label{carr\'es}
\end{figure}
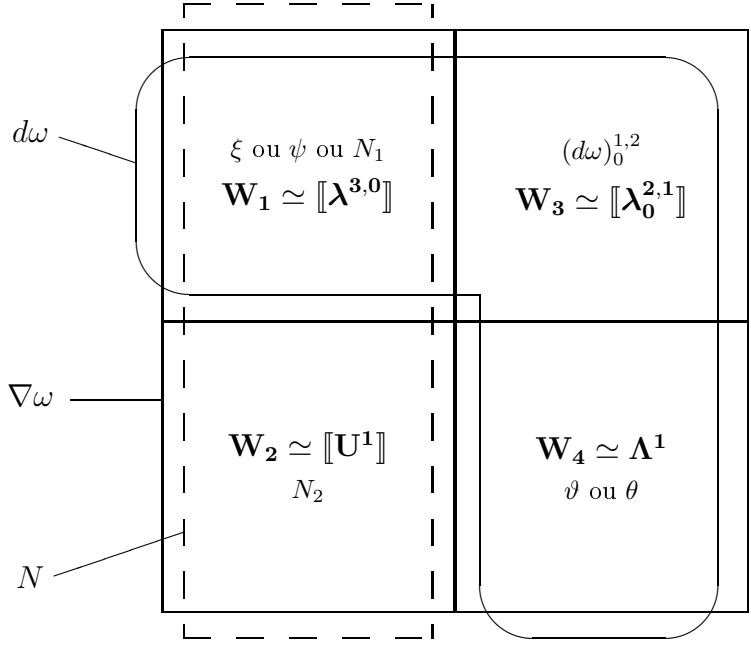

La figure \ref{carr\'es} r\'esume les identifications qu'on peut faire entre les composantes de $\n \omega$, $d\omega$, $\overline \delta$ et du tenseur de Nijenhuis $N$ (voir plus loin). Plus de d\'etails (notamment les isomorphismes $U(n)$-invariants qui sous-tendent ces identifications) sont donn\'es au chapitre 3 de la th\`ese \cite{bu3}.

\vs

Les vari\'et\'es de type $W_1 \oplus W_2 \oplus W_4$
sont caract\'eris\'ees par
\[ d\omega = \psi + \omega \wedge \theta \]
et parmi celles-ci, les vari\'et\'es de type $W_1 \oplus W_4$
sont celles ayant
\begin{equation}
\overline \delta = \xi + \vartheta
\label{xi+theta'}
\end{equation}

A titre de comparaison, les vari\'et\'es de type $W_2 \oplus
W_4$ sont caract\'eris\'ees par
\[ d\omega = \omega \wedge \theta \]
De cette derni\`ere \'equation on tire, en diff\'erentiant
\[ 0 = d^2\omega = \omega \wedge d\theta \]
De l\`a que, en dimension sup\'erieure \`a 6, $\theta$ est
ferm\'ee (localement exacte) et les vari\'et\'es de la classe $W_2
\oplus W_4$ sont localement conformes \`a des vari\'et\'es de type
$W_2$. En particulier $\n \omega \in {\bf W}_4$ si et seulement si la
vari\'et\'e est conform\'ement k{\"a}hlerienne, comme il est bien connu. En dimension 4 a contrario, {\it toutes} les vari\'et\'es presque hermitiennes $M$ sont de type $W_2 \oplus W_4$ et $M$ appartient \`a $W_4$ si et seulement si la structure presque complexe est int\'egrable.

On donne maintenant la d\'efinition suivante :

\begin{defi}
On appelle \emph{sp\'eciale} une vari\'et\'e presque hermitienne
telle que la norme de $\psi=(d\omega)^{3,0}$ est constante, \'egale \`a 1.
\end{defi}

On peut toujours supposer, pour notre probl\`eme, qu'une
vari\'et\'e de type $W_1 \oplus W_4$ est sp\'eciale, localement sur un ouvert o\`u $\psi$ ne s'annule pas, en effectuant un changement conforme appropri\'e. En effet, la classe $W_1 \oplus W_4$ est stable par changement conforme (voir \cite{gr4}). Soit $f$ une fonction strictement positive d\'efinie sur un ouvert de $M$. Si on substitue $g'=fg$ \`a $g$, la forme de K\"ahler $\omega$ est remplac\'ee par $\omega'=f\omega$ telle que
\beq
d\omega' & = & df \wedge \omega + fd\omega \\
         & = & (\frac{df}{f}+\theta) \wedge \omega' + f\psi
\eeq
Par cons\'equent, $\theta$ est remplac\'ee par la 1-forme cohomologue
$\theta'=\theta + \frac{df}{f}$ (en particulier $\theta$
est ferm\'ee si et seulement si $\theta'$ l'est), $\psi$ par $\psi'=f\psi$ et la norme de $\psi$ par
\[
| \psi' |_{g'}=f| \psi |_{g'}= f^{-\frac{1}{2}} | \psi |_g
\]
On peut donc choisir $f$, sur un ouvert o\`u $\psi$ ne s'annule pas, pour
que la norme de $\psi'=(d\omega')^{3,0}$ par rapport \`a la nouvelle
m\'etrique soit constante \'egale \`a 1.

En fait on peut le supposer globalement. Plus pr\'ecis\'ement :
\begin{lemm}
Les propositions suivantes sont \'equivalentes : \\
(\romannumeral 1) Les vari\'et\'es de type $W_1 \oplus W_4$, de
dimension 6, sont localement conform\'ement NK. \\
(\romannumeral 2) Les vari\'et\'es de type $W_1 \oplus W_4$
\emph{sp\'eciales} de dimension 6 sont NK.
\label{psi=1}
\end{lemm}

\begin{proof}
D'apr\`es ce qui pr\'ec\`ede, le changement conforme est d\'etermin\'e par le rapport des normes de $\psi$, $\psi'$. Or les vari\'et\'es NK non k\"ahleriennes en dimension 6 sont de type constant, d'apr\`es \cite{gr3}, c'est-\`a-dire ont $\n \omega$ ou $d\omega$ (qui dans ce cas est de type pur
(3,0)+(0,3), \'egale \`a $\psi$) de norme constante non nulle. D\`es lors, si une vari\'et\'e de type $W_1 \oplus W_4$ est localement conforme \`a une vari\'et\'e NK, le changement conforme est donn\'e par un multiple de $f = | \psi |^2$. En particulier une vari\'et\'e sp\'eciale de type $W_1 \oplus W_4$ est localement conform\'ement NK si et seulement si elle est NK.

Montrons $(\romannumeral 2) \Rightarrow (\romannumeral 1)$. Soit $U$ l'ouvert de $M$ o\`u $\psi$ ne s'annule pas. La structure presque hermitienne sur $U$ est conforme \`a une structure presque hermitienne sp\'eciale de type $W_1 \oplus W_4$, c'est-\`a-dire, par (\romannumeral 2), NK. Autrement dit la forme de K\"ahler est exacte ($\theta = \frac{df}{f}$) donc ferm\'ee sur $U$. Mais d'autre part, sur le compl\'ementaire de $U$, $\psi=0$ donc la structure presque hermitienne est de type $W_4$, $d\omega = \theta \wedge \omega$. Cela implique comme pr\'ec\'edemment $d\theta = 0$ sur tout ouvert de $M \setminus U$. Finalement $d\theta$ est nulle sur tout $M$ par continuit\'e. 
\end{proof}

\begin{NB}
Ce dernier argument est pr\'esent\'e de fa\c con \`a peine diff\'erente dans l'article de Cleyton, Ivanov \cite{cl4}. Les auteurs montrent en outre (lemma 8) qu'une vari\'et\'e $W_1 \oplus W_4$ (non localement conform\'ement k\"ahlerienne) est localement conforme \`a une vari\'et\'e NK si et seulement si elle est \emph{globalement} conforme \`a une vari\'et\'e NK. Ajout\'e au th\'eor\`eme \ref{W1+W4=W1}, cela signifie qu'il n'y a \emph{aucune} diff\'erence entre la classe des vari\'et\'es strictement NK et $W_1 \oplus W_4 \setminus W_4$ \`a une fonction pr\`es $f$ d\'efinie sur tout $M$. Les m\^emes r\'esultats sont valables, par \cite{ma2,cl4}, en rempla\c cant les vari\'et\'es \og nearly K\"ahler \fg \ par les vari\'et\'es presque parall\`eles $G_2$.
\end{NB}

On se restreint dans la suite aux vari\'et\'es de type $W_1 \oplus
W_4$ sp\'eciales de dimension 6. L'int\'er\^et est qu'on peut
s'appuyer sur cette 3-forme $\psi$ qui
ne s'annule jamais pour obtenir des restrictions sur la g\'eom\'etrie
de la vari\'et\'e, en s'inspirant des m\'ethodes utilis\'ees pour les
vari\'et\'es strictement NK. En premier lieu on
dispose maintenant sur la vari\'et\'e d'une structure $SU(3)$.

Rappelons quelques propri\'et\'es
alg\'ebriques des formes de type (3,0)+(0,3) en dimension 6. Quelle que soit la 3-forme $\psi$ de type (3,0)+(0,3), il existe une unique $\phi \in \leftr \llambda^{3,0} \rightr$ telle que
\[ \Psi = \psi +i\phi \]
est de type (3,0). Si $\psi$ est de norme constante non nulle, $\Psi$ est
une 3-forme volume complexe, d\'efinissant une r\'eduction de $M$ \`a $SU(3)$, et $\psi$, $\phi$ forment une base orthogonale de $\leftr \llambda^{3,0} \rightr$ pour le produit scalaire induit par la m\'etrique.

\begin{lemm} \ \\
(\romannumeral 1) $\forall X \in TM, \quad \iota_X\psi =
    \iota_{JX}\phi$ \\
(\romannumeral 2)
$\omega \wedge \iota_X \psi = -JX^{\flat} \wedge \psi, \ $
et de m\^eme pour $\phi$. \\
(\romannumeral 3)
$\forall \eta \in \LLambda^1, \quad \eta \wedge \psi = J\eta \wedge \phi$
\label{psi,phi}
\end{lemm}

Une autre fa\c con de d\'efinir $\phi$ \`a partir de $\psi$ est
\begin{equation}
\phi = * \psi
\label{*psi}
\end{equation}

Pour une vari\'et\'e presque hermitienne sp\'eciale, on a donc une structure $SU(3)$ naturelle, compatible avec la structure $U(3)$, associ\'ee \`a $\psi = (d\omega)^{3,0}$. Or dans ce cas, F. M. Cabrera a remarqu\'e dans \cite{ma} -- et Salamon, Chiossi \cite{ch} ont montr\'e que toute l'information pour la torsion intrins\`eque est contenue dans les diff\'erentielles de $\omega$, $\psi$, $\phi$. En particulier, les composantes de $\overline \delta$ ou $\nabla \omega$ dans ${\bf W}_1 \oplus {\bf W}_2$ sont d\'etermin\'ees par la partie de type (2,2) de $d\psi$, $d\phi \in \LLambda^4$. Plus pr\'ecis\'ement, $[\lambda^{2,2}_0]$ est r\'eduit \`a \{0\} en dimension 6, par cons\'equent $[\lambda^{2,2}] \simeq [\lambda^{1,1}] \simeq [\lambda^{1,1}_0]
\oplus \RM$ et il existe $n_1$, $m_1 \in
\RM$ et $\nu_2$, $\mu_2 \in [\llambda^{1,1}_0]$ telles que
\[ (d\psi)^{2,2} = n_1 \omega \wedge \omega + \nu_2 \wedge \omega \]
\[ (d\phi)^{2,2} = m_1 \omega \wedge \omega + \mu_2 \wedge \omega \]

Alors,
\begin{prop}[Salamon, Chiossi]
Soit $M$ une vari\'et\'e riemannienne munie d'une structure $SU(3)$
d\'efinie par $(\omega,\psi,\phi)$. \\
1) La composante de la torsion intrins\`eque dans ${\bf W}_1$ est totalement d\'etermin\'ee par $n_1$, $m_1$ c'est-\`a-dire $w_1=0$ si et seulement si $n_1 = m_1 = 0$. \\
2) De m\^eme la composante de $\overline \delta$ ou $\nabla \omega$ suivant ${\bf W}_2$ est totalement d\'etermin\'ee par $\nu_2$, $\mu_2$.
\label{n2,m2}
\end{prop}

\begin{proof}
Cela r\'esulte premi\`erement de raisons th\'eoriques concernant les
repr\'esentations de $SU(3)$. L'application qui au 1-jet de la
structure $U(3)$ (repr\'esent\'e par $\n \omega$ ou la torsion
intrins\`eque $\overline \delta$) associe le quadruplet $(w_1,w_2,w_3,w_4)$
est une bijection. Les sous-espaces $W_1$ et $W_2$ se
d\'ecomposent sous l'action de 
$SU(3)$ en
\[ W_1 = \RM \oplus \RM \]
\[ W_2 = \mathfrak{su}(3) \oplus \mathfrak{su}(3) \]
Au contraire $W_3$, $W_4$ sont toujours irr\'eductibles (et aucun
d'eux n'est isomorphe \`a $\RM$ ni $\mathfrak{su}(3)$). On
voit $n_1$, $m_1$ (resp. $\nu_2$, $\mu_2$) comme des applications de
l'espace des 1-jets de structures $SU(3)$ (faisant intervenir la
d\'ecomposition en types de $d\psi$, $d\phi$) dans $\RM$
(resp. $[\lambda^{1,1}_0] \simeq \mathfrak{su}(3)$). Maintenant
l'espace des 1-jets de structures
$SU(3)$ est envoy\'e surjectivement sur l'espace des 1-jets de
structures $U(3)$ et par le lemme de Schur, toutes ces applications
\'etant $SU(3)$-invariantes, on obtient le r\'esultat
annonc\'e.

Une autre fa\c con de le voir est la suivante. Soit $\alpha$ une forme diff\'erentielle de type
$(p,q)$. La partie de type
$(p-1,q+2)$ de $d\alpha$ est donn\'ee par le tenseur qui mesure le d\'efaut
d'int\'egrabilit\'e de la structure presque complexe c'est-\`a-dire $N$, le tenseur
de Nijenhuis. Concr\`etement si $\alpha \in \llambda^{1,0}$
\[
(d\alpha)^{0,2}(X,Y)=d\alpha(X^{0,1},Y^{0,1})
=-\alpha([X^{0,1},Y^{0,1}])
\]
D'o\`u, par (\ref{N}) :
\[
(d\alpha)^{0,2}(X,Y)=-\alpha(N(X,Y))
\]
Par l\`a le tenseur $N$ peut-\^etre vu comme un \'el\'ement de
$Hom(\llambda^{1,0},\llambda^{0,2}) \simeq \llambda^{0,1} \otimes
\llambda^{0,2}$, ou autrement dit le tenseur {\it r\'eel} $N$ vit dans un
espace isomorphe \`a $\leftr \llambda^{0,1} \otimes
\llambda^{0,2} \rightr \simeq {\bf W}_1 \oplus {\bf W}_2$. On le d\'ecompose en cons\'equence
en
\begin{equation}
N = N_1 + N_2 \label{N12}
\end{equation}
On \'etablit de m\^eme une formule plus g\'en\'erale, pour toute
forme de type $(p,0)$ mais aussi $(p,0)+(0,p)$ :
\begin{equation}
(d\alpha)^{p-1,2}= N \# \alpha
\label{Ndiese}
\end{equation}
o\`u $\# : (\LLambda^2 \otimes TM) \times
\LLambda^p \to \LLambda^{p+1}$. Maintenant $N_1$ est li\'e \`a $w_1$ ou $\psi$ et on peut montrer, quelle que soit $\alpha \in \leftr \llambda^{3,0} \rightr$ :
\begin{eqnarray}
(d\alpha)^{2,2} & = & N_1 \# \alpha + N_2 \# \alpha \nonumber \\
                & = & -\frac{2}{3} \la \phi,\alpha \ra \omega \wedge
                \omega + \gamma \wedge \omega \label{N*a}
\end{eqnarray}
o\`u $\gamma \in [\lambda^{1,1}_0]$. En prenant $\alpha$ \'egale
\`a $\psi$ ou
$\phi$ on obtient ainsi, puisqu'elles forment une base de $\leftr
\llambda^{3,0} \rightr$ que $N_1$ est enti\`erement d\'etermin\'e par
$n_1$, $m_1$ et $N_2$ est d\'etermin\'e par $\nu_2$, $\mu_2$.
\end{proof}

On a m\^eme un r\'esultat plus pr\'ecis du fait que la structure
$SU(3)$ n'est pas quelconque mais que la 3-forme $\psi$ est issue de
$d\omega$ comme sa partie de type (3,0)+(0,3).

\begin{lemm}
Pour une vari\'et\'e presque hermitienne sp\'eciale de dimension
6, $n_1=0$ c'est-\`a-dire
\[ (d\psi)^{2,2}=\nu_2 \wedge \omega, \]
o\`u $\nu_2 \in [\llambda^{1,1}_0]$ et $m_1$ est constante.
\label{n1}
\end{lemm}

\begin{proof}
On a directement, par (\ref{N*a}), $n_1 = -\frac{2}{3} \la \phi,\psi \ra = 0$ et
$m_1 = -\frac{2}{3} | \phi | ^2 = -\frac{2}{3}$.
\end{proof}

\begin{prop}
Soit $M$ une vari\'et\'e presque hermitienne de type $W_1 \oplus W_4$
sp\'eciale de dimension 6. Il existe $\sigma \in \LLambda^1$ telle que
\begin{equation}
d\psi = \sigma \wedge \psi
\label{dpsi}
\end{equation}
\begin{equation}
d\phi = \sigma \wedge \phi - \frac{2}{3}\omega \wedge \omega
\label{dphi}
\end{equation}
\label{dpsi,dphi}
\end{prop}

\begin{proof}
D'abord par la proposition \ref{n2,m2} puis le lemme \ref{n1}, pour
une vari\'et\'e presque hermitienne sp\'eciale, $w_2=0$ si et
seulement si il existe $\sigma$ et $\varsigma \in \LLambda^1$ telles que
\[ d\psi = \sigma \wedge \psi \]
\[ d\phi = \varsigma \wedge \phi - \frac{2}{3} \omega \wedge \omega \]
En effet en dimension 6, $\leftr \llambda^{3,0} \rightr$ est de
dimension 2 mais gr\^ace \`a l'\'equation (\romannumeral 2) du
lemme \ref{psi,phi} toute forme de type (3,1)+(1,3) est
d\'ecomposable en le produit ext\'rieur d'une parmi $\psi$, $\phi$
avec une certaine 1-forme. Cela signifie que les applications

\[
\begin{array}{ccccccc}
\LLambda^1 & \to & \leftr \llambda^{3,1} \rightr
& \text{ et } & \LLambda^1 & \to & \leftr \llambda^{3,1} \rightr \\
\eta & \mapsto & \eta \wedge \psi & & \eta & \mapsto & \eta \wedge
\psi
\end{array}
\]

sont des isomorphismes.

En outre ces 1-formes ne sont pas arbitraires : on peut calculer $\sigma = \varsigma$ (voir \cite{bu3}, proposition 3.3.16). Il s'agit en r\'ealit\'e d'un fait plus g\'en\'eral, qui concerne toutes les vari\'et\'es $SU(3)$ d'apr\`es \cite{ch} : il existe une 1-forme $\sigma$ telle que $(d\psi)^{3,1}=\sigma \wedge \psi$ et $(d\phi)^{3,1}=\sigma \wedge \phi$.
\end{proof}

On peut maintenant \'enoncer le
\begin{theo}
Soit $M$ une vari\'et\'e de type $W_1 \oplus W_2 \oplus W_4$
c'est-\`a-dire v\'erifiant
\begin{equation}
d\omega = \psi + \omega \wedge \theta
\label{W124}
\end{equation}
Alors $(d\psi)^{2,2}=0$ si et seulement si
$(d\theta)^{1,1}=0$ et si $M$ est de dimension 6, $w_2=0$,
i.e. la vari\'et\'e est de type $W_1 \oplus W_4$ si et seulement si elle est
localement conform\'ement NK.
\label{W1=W1+W4sp}
\end{theo}

\begin{proof}
Le premier point r\'esulte simplement de la diff\'erentiation de
l'\'equation caract\'eristique (\ref{W124})
\begin{equation}
0 = d\psi + \psi \wedge \theta + \omega \wedge d\theta
\label{ddomega}
\end{equation}
La 4-forme $\omega \wedge d\theta$ se scinde suivant les types en
\[ \omega \wedge d\theta = \omega \wedge (d\theta)^{2,0} + \omega
\wedge (d\theta)^{1,1}
\]
Le dernier terme doit \^etre nul si $d\psi$ est de type pur
(3,1)+(1,3).

Pour prouver le deuxi\`eme point on se ram\`ene au cas des
vari\'et\'es sp\'eciales gr\^ace au lemme \ref{psi=1}. Si
$w_2=0$, par les lemmes \ref{n2,m2} et \ref{n1}, $d\psi$ est de
type (3,1)+(1,3) d'o\`u $(d\theta)^{1,1}=0$ par ce qui pr\'ec\`ede
et on se trouve dans la situation de la proposition \ref{dpsi,dphi}. En diff\'erentiant (\ref{dpsi}) on obtient
premi\`erement
\[ 0 = d^2\psi = d\sigma \wedge \psi \]
ce qui implique que $d\sigma$ est de type (1,1) car l'application
\beq
\leftr \llambda^{2,0} \rightr & \to     & \LLambda^5 \\
\eta                         & \mapsto & \eta \wedge \psi \eeq est
bijective (NB : en fin de compte, on a obtenu tous les isomorphismes de repr\'esentations de $SU(3)$ suivants : $\leftr \lambda^{3,1} \rightr
\simeq \Lambda^1 \simeq \leftr \lambda^{2,0} \rightr \simeq
\Lambda^5$.) Mais alors on a aussi $d\sigma \wedge \phi = 0$ (il
n'y a pas de forme de type (4,1)) et en diff\'erentiant
(\ref{dphi}), en tenant compte de $\psi \wedge \omega =0$,

\beq
0 & = & -\sigma \wedge d\phi - \frac{4}{3}d\omega \wedge \omega \\
  & = & \frac{2}{3}\sigma \wedge \omega \wedge \omega - \frac{4}{3} \theta
  \wedge \omega \wedge \omega \\
  & = & \frac{2}{3}(\sigma - 2\theta) \wedge \omega \wedge \omega
\eeq
Par cons\'equent $\theta = \frac{1}{2}\sigma$ et $d\theta$ est de
type (1,1) c'est-\`a-dire, finalement, nulle.

Comme on a suppos\'e que la vari\'et\'e \'etait
sp\'eciale, on doit montrer que cela implique que
$\theta$ aussi est nulle. Or, (\ref{ddomega}) devient maintenant
$d\psi = \theta \wedge \psi$ et (\ref{dpsi}), $d\psi = 2\theta
\wedge \psi$ d'o\`u $\theta = 0$.

R\'eciproquement Gray et Hervella \cite{gr4} ont montr\'e qu'une
vari\'et\'e localement conform\'ement NK est de type $W_1 \oplus W_4$.

Les \'equations (\ref{dpsi}) et (\ref{dphi}) sont finalement les \'equations
caract\'eristiques d'une vari\'et\'e
strictement NK de dimension 6 obtenues par Reyes-Carrion \cite{re} (voir
aussi \cite{hi3}). Toute vari\'et\'e presque hermitienne de type $W_1 \oplus
W_4$ sp\'eciale de dimension 6 est strictement NK et par le lemme
\ref{psi=1} toute vari\'et\'e de type $W_1 \oplus W_4$ de dimension 6
est localement conforme \`a une vari\'et\'e NK.
\end{proof}

De ce th\'eor\`eme et du corollaire \ref{S6} d\'ecoule le th\'eor\`eme
annonc\'e \ref{twisteurs}. Mais c'est aussi une pi\`ece d'un
probl\`eme plus g\'en\'eral.

\section{Vari\'et\'es presque hermitiennes conformes}

Parmi les 16 classes de vari\'et\'es presque hermitiennes $\bigoplus_{i
\in I} W_i$, $I \subset \{1,2,3,4\}$ qu'ils d\'efinissent (voir d\'efinition \ref{Gray Hervella}), Gray et Hervella \cite{gr4} ont remarqu\'e que 8 -- celles dont la
d\'efinition comprend $W_4$ -- sont invariantes conformes. En particulier, $\forall I \subset \{1,2,3\}$, $\bigoplus_{i \in I} W_i \oplus W_4$ contient non seulement les vari\'et\'es de de la classe plus
petite $\bigoplus_{i \in I} W_i$ mais encore toutes les
vari\'et\'es localement conformes \`a des vari\'et\'es de ce type.
La question se pose alors si ces derni\`eres \'epuisent toute la
classe. 

Par le th\'eor\`eme \ref{W1=W1+W4sp} la r\'eponse est positive
en dimension 6
dans le cas de $W_1 \oplus W_4$. On a vu au d\'ebut de la section 5
que $W_2 \oplus W_4$ \'etait
de m\^eme constitu\'ee uniquement de vari\'et\'es localement
conformes \`a des vari\'et\'es $W_2$, c'est-\`a-dire presque k\"ahleriennes ou symplectiques, et
que $W_4$ ne contenait que les vari\'et\'es localement
conform\'ement k\"ahleriennes en dimension sup\'erieure \`a 6.

Au
contraire il existe des vari\'et\'es \emph{hermitiennes}, c'est-\`a-dire de type
$W_3 \oplus W_4$, en dimension 6, non localement conformes \`a des vari\'et\'es $W_3$. En effet soit $M$ le produit d'une surface de Riemann $S$,
de forme de K\"ahler $\omega_0$, et d'une vari\'et\'e
hermitienne $N$ de dimension 4 (une surface complexe) non LCK (c'est
le cas g\'en\'eriquement.) La forme de K\"ahler $\omega_1$ de $N$
v\'erifie $d\omega_1 = \theta_1 \wedge \omega_1$, pour une 1-forme
$\theta_1$ non ferm\'ee par hypoth\`ese. Alors la forme de K\"ahler de $M =
S \times N$ est $\omega = \omega_0 + \omega_1$ et on v\'erifie que la
forme de Lee est $\theta = \frac{1}{2} \theta_1$, de sorte que $d\theta
\neq 0$.

De plus on s'appr\^ete \`a donner (voir corollaire \ref{W124>W12}) un r\'esultat d'existence locale de
vari\'et\'es de type $W_1 \oplus
W_2 \oplus W_4$ de dimension 6 non issues de vari\'et\'es
semi-k\"ahleriennes c'est-\`a-dire
v\'erifiant (\ref{W124}) mais cette fois encore $d\theta \neq 0$.

Afin de donner
une formulation intrins\`eque de ces r\'esultats, on consid\`ere des
vari\'et\'es presque hermitiennes conformes, au sens suivant. Soit $M$
une vari\'et\'e de dimension $m$. On note
$\mathcal L^k$, pour tout entier $k$, et on appelle \emph{fibr\'e des scalaires de poids $k$},
le fibr\'e associ\'e de $GL(M)$ :
\[ \mathcal L^k = GL(M) \times_{|det|^{k/n}} \RM, \]
de sorte que $\mathcal L^0$ s'identifie \`a $M
\times \RM$ et $\mathcal L^k \otimes \mathcal L^l = \mathcal L^{k+l}$.
Une structure conforme sur $M$ est la donn\'ee d'une classe conforme
de m\'etriques $C$ ou bien encore d'une section $c$ de $S^2(T^*M)
\otimes \mathcal L^{-2}$ telle que pour $X \in TM$, non nul,
$c(X,X)$ est strictement positif (i.e. s'\'ecrit $l \otimes
l$ pour une section non nulle $l$ de $\mathcal L$). Le lien entre les deux
d\'efinitions est le suivant : tout choix d'une m\'etrique $g$ dans $C$ correspond
\`a une trivialisation ou une section $l$ de $\mathcal L$ qui permette
d'\'ecrire
\[ g = c \otimes l^2
\]

La structure
conforme $c$ fournit les isomorphismes
\[ 
\begin{array}{ccccccc}
TM & \to & \LLambda^1 \otimes \mathcal L^2 & \quad \text{et} \quad & \LLambda^1 & \to & TM \otimes \mathcal L^{-2} \\
X  & \mapsto & X^{c} & & \eta & \mapsto & \eta^{\diamond}
\end{array}
\]
analogues des \emph{isomorphismes musicaux} $X \mapsto X^{\flat}$ et $\eta
\mapsto \eta^{\sharp}$, dans le cas riemannien. La compatibilit\'e
entre une m\'etrique et une structure presque complexe, telle
qu'elle d\'efinit une vari\'et\'e presque hermitienne, si elle a
lieu pour une m\'etrique de $C$, a lieu pour toutes les
m\'etriques de la classe. C'est pourquoi on donne la d\'efinition suivante :
\begin{defi}
On appelle vari\'et\'e presque hermitienne conforme une
vari\'et\'e $M$ munie d'une structure conforme $c$ et d'une structure
presque complexe $J$ compatibles c'est-\`a-dire satisfaisant
\[ \forall X,Y \in TM, \quad c(JX,JY)=c(X,Y) \]
\end{defi}
De fa\c con \'equivalente il existe une r\'eduction du fibr\'e des
rep\`eres \`a $CU(n)$, identifi\'e canoniquement \`a $\RM^*_+ \times
U(n)$. 

En l'absence de m\'etrique,
on ne dispose plus de connexion de Levi-Civita mais uniquement de connexions de
Weyl. Soit $D$ la d\'eriv\'ee covariante d'une telle connexion :
lin\'eaire, sans torsion et pr\'eservant $c$. On consid\`ere $DJ \in
\LLambda^1 \otimes \mathfrak{co}(M)$. En fait, $\forall X \in TM$, $D_X
J$ est antisym\'etrique (relativement \`a $c$ ou \`a toute m\'etrique de
la classe conforme), c'est-\`a-dire appartient \`a $\mathfrak{so}(M)$, bien
d\'efini m\^eme en l'absence de m\'etrique. De plus, en d\'erivant $J^2
= -Id$, on obtient qu'il
anticommute \`a $J$. Finalement on a comme avant $DJ \in \LLambda^1 \otimes \mathfrak{u}(M)^{\perp}$, o\`u
l'orthogonal est pris dans $\mathfrak{so}(M)$. En
restreignant la repr\'esentation de
$CU(n)$ \`a $U(n)$, on d\'ecompose
ce fibr\'e en composantes irr\'eductibles en chaque point comme en (\ref{W1234}).
Maintenant, quel que soit $i = 1,2,3,4$, $W_i$ est invariant sous l'action de $\RM^*_+$
par cons\'equent il s'agit de la d\'ecomposition irr\'eductible de la repr\'esentation de $CU(n)$. On souhaite obtenir une classification qui ne d\'epende pas de la connexion $D$. Soit $D'$ une connexion de Weyl diff\'erente de $D$. On cherche \`a comparer $DJ$ et $D'J$. Chaque connexion $D$, $D'$ induit une connexion lin\'eaire sur $\mathcal L$, not\'ee $\mathcal D$, $\mathcal D'$. R\'eciproquement, on peut calculer $D$, $D'$ \`a partir de $\mathcal D$, $\mathcal D'$ par un analogue de la
formule de Koszul :
\begin{eqnarray*}
2c(D_X Y Z) & = & \mathcal D_X(c(Y,Z))+\mathcal D_Y(c(X,Z))-\mathcal
D_Z(c(X,Y)) \nonumber \\
& & + c(Z,[X,Y]) - c(Y,[X,Z]) - c(X,[Y,Z])
\end{eqnarray*}
On appelle $\tau$ la 1-forme telle que, pour $l \in \mathcal L$
\[ (\mathcal D_X - {\mathcal D'}_X)l = 2\tau(X)l \]
Alors
\begin{equation}
(D-D')_X Y = \tau(X)Y + \tau(Y)X - c(X,Y)\tau^{\diamond}
\label{D-D'}
\end{equation}
Puis en identifiant $\mathfrak{so}(n)$ et $\LLambda^2 \otimes \mathcal L^2$,
\begin{equation}
(D-D')_X J = X^c \wedge J\tau + JX^c \wedge \tau
\label{D-D'J}
\end{equation}
Par cons\'equent les trois premi\`eres composantes de $DJ$ ne changent pas, lorsqu'on modifie la connexion de Weyl, mais seulement la quatri\`eme, \`a valeurs dans ${\bf W}_4$. C'est exactement
ce qu'ont d\'emontr\'e Gray et Hervella en se restreignant aux structures de Weyl ferm\'ees, c'est-\`a-dire co\"incidant localement avec la connexion de Levi-Civita d'une m\'etrique de $C$. D'ailleurs on peut
toujours modifier $D$ pour que la composante suivant ${\bf W}_4$ soit nulle :
\begin{propdefi}
Soit $(M,c,J)$ une vari\'et\'e presque hermitienne con\-forme. Il
existe une unique connexion de Weyl, not\'ee $D^J$, telle que la
composante de $D^J J$ dans ${\bf W}_4$ est nulle. On appelle cette
connexion la connexion adapt\'ee \`a $J$.
\end{propdefi}
\begin{proof}
En effet pour une connexion de Weyl quelconque $D$, il existe une 1-forme
$\eta$ telle que la composante dans ${\bf W}_4$ de $DJ$ s'\'ecrit (voir (\ref{(theta)}) pour l'analogue riemannien)
\[ L(D) = JX^{\flat} \wedge \eta + X^{\flat} \wedge J\eta
\]
Il suffit alors de poser
\[ D^J_X Y = D_X Y - \eta(X)Y - \eta(Y)Z + c(X,Y)\eta^{\diamond} \]
en s'inspirant de (\ref{D-D'}).
\end{proof}

La d\'efinition suivante ne d\'epend pas du choix de $D$ :
\begin{defi}
Une vari\'et\'e presque hermitienne conforme $(M,c,J)$ est dite de type $W^c$,
o\`u $W=\bigoplus_{i \in I} W_i$, $I \subset \{1,2,3\}$, si pour
toute connexion de Weyl $D$, $DJ$ est une section de ${\bf W} \oplus {\bf W}_4$.
\end{defi}

Bien s\^ur on a
\begin{prop}
Soient $(M,c,J)$ une vari\'et\'e presque hermitienne conforme, $g \in C$,
la classe conforme d\'efinie par $c$. La vari\'et\'e presque
hermitienne $(M,g,J)$ est de type $\bigoplus_{i \in I} W_i \oplus W_4$
si et seulement si $(M,c,J)$ est de type $\bigoplus_{i \in I} W_i^c$.
\end{prop}

De plus,
\begin{defi}
Une vari\'et\'e presque hermitienne conforme $(M,c,J)$ est dite
ferm\'ee (resp. exacte) si la connexion adapt\'ee \`a $J$
est ferm\'ee (resp. exacte) en tant que structure de Weyl.
\end{defi}

Rappelons qu'une structure de Weyl $D$ est dite ferm\'ee (resp. exacte) si
pour une m\'etrique $g \in C$, la 1-forme $\eta$, appel\'ee forme de Lee de $(D,g)$, qui
mesure la diff\'erence avec la connexion de Levi-Civita comme en (\ref{D-D'}) est ferm\'ee (resp. exacte). Cette d\'efinition ne d\'epend pas du choix de la m\'etrique dans $C$, en effet la forme de Lee de
$(D,g'=e^{2f}g)$ est \'egale \`a $\eta'=\eta + df$, cohomologue \`a $\eta$. Maintenant, pour la connexion adapt\'ee \`a $J$, la forme de Lee de
$(D^J,g)$ est \'egale \`a $\frac{1}{2}\theta$ c'est-\`a-dire (\`a un
facteur pr\`es) \`a la forme de Lee de la vari\'et\'e presque
hermitienne $(M,g,J)$. En effet par (\ref{(theta)})
\[ (\n - D^J)J = L(\n) = \frac{1}{2}(X^c \wedge J\theta + JX^c \wedge
\theta)
\]
On peut alors reformuler le th\'eor\`eme \ref{W1+W4=W1}

\begin{theo}
Toute vari\'et\'e presque hermitienne conforme de dimension 6, de type $W_1^c$
est ferm\'ee.
\end{theo}

Par une remarque pr\'ec\'edente, on a le m\^eme r\'esultat pour les vari\'et\'es de la classe $W_2^c$, en dimension sup\'erieure \`a
6. Dans la section suivante on consid\`ere le cas des vari\'et\'es
presque hermitiennes conformes de type $W_1^c \oplus W_2^c$.

\section{Vari\'et\'es de type $W_1 \oplus W_2 \oplus W_4$} 

Soit $V$ un espace vectoriel de dimension 6. Tout part de
l'observation faite par Hitchin que la donn\'ee de deux formes $\omega
\in \Lambda^2 V^*$ et $\psi \in \Lambda^3 V^*$ est suffisante pour d\'efinir
une action de $SU(3)$ sur $V$ (ou au niveau des vari\'et\'es, la
donn\'ee de deux formes {\it diff\'erentielles} $\omega$ et $\psi$ suffit
\`a d\'efinir une {\it r\'eduction} du fibr\'e des rep\`eres \`a $SU(3)$) moyennant certaines
conditions de r\'egularit\'e -- pour chaque forme s\'epar\'ement -- et
de compatibilit\'e. C'est-\`a-dire que dans cette donn\'ee et sous ces
conditions est en fait incluse la donn\'ee de l'endomorphisme $J$ de
carr\'e -1 (ou de la structure presque complexe) et du produit
scalaire (ou de la m\'etrique). Voici comment on s'y prend.

Tout d'abord $\psi$ doit \^etre une 3-forme \emph{stable}, plus pr\'ecis\'ement on demande que son stabilisateur soit isomorphe
\`a $SL(3,\CM)$ (induisant un isomorphisme $V \simeq \CM^3$). Pour caract\'eriser les 3-formes stables, Hitchin \cite{hi3}
d\'efinit une certaine quantit\'e $\kappa : \Lambda^3 V^* \to
(\Lambda^6 V^*)^2$ qui doit v\'erifier dans notre cas :
\begin{equation}
\kappa(\psi) < 0
\tag{r1}
\label{r1}
\end{equation}

La 3-forme $\psi$ d\'efinit alors une structure presque complexe $J$ par rapport \`a laquelle elle est de type
  $(3,0)+(0,3)$. D\`es lors la premi\`ere condition de compatibilit\'e
  est que $\omega$ soit de type (1,1), ce qui \'equivaut \`a :
\begin{equation}
\omega \wedge \psi = 0
\tag{c1}
\label{c1}
\end{equation}

Comme pr\'ec\'edemment on associe \`a $\psi$ une 3-forme $\phi$ de fa\c con unique pour
que $\psi + i\phi$ soit de type (3,0). Une cons\'equence de (\ref{r1}) est que $\psi$ est \emph{non-d\'eg\'en\'er\'ee} c'est-\`a-dire que l'\'el\'ement de volume canoniquement associ\'e \`a
$\psi$, $\mu(\psi) = \psi \wedge \phi \in \Lambda^6 V^*$ est non nul. On demande aussi que
$\omega$ soit non-d\'eg\'en\'er\'ee :
\begin{equation}
\mu(\omega) = \omega \wedge \omega \wedge \omega \neq 0
\tag{r2}
\label{r2}
\end{equation}

La deuxi\`eme condition de compatibilit\'e est que la forme
bilin\'eaire sym\'etrique d\'efinie par $J$ et $\omega$ soit d\'efinie
positive :
\begin{equation}
(X,Y) \mapsto g(X,Y) =
\omega(X,JY) \ > 0
\tag{c2}
\label{c2}
\end{equation}

Enfin la derni\`ere est que $\psi$ soit de norme 1 pour $g$, ce qui
dans le langage des formes s'\'ecrit :
\begin{equation}
\mu(\psi) = \frac{2}{3} \mu(\omega)
\tag{c3}
\label{c3}
\end{equation}

\begin{prop}
L'ensemble des couples $(\omega,\psi)$ v\'erifiant les conditions de
r\'egularit\'e (\ref{r1}) et (\ref{r2}) et la condition de
compatibilit\'e (\ref{c2}) est un ouvert non vide de $\Lambda^2 V^*
\times \Lambda^3 V^*$
\label{open}
\end{prop}

\begin{proof}
Cela r\'esulte de la d\'efinition choisie par Hitchin d'une $p$-forme stable : l'orbite de $\psi$ sous $GL(V)$ est ouverte dans $\Lambda^p V^*$. En dimension 6, $\Lambda^3 V^*$ contient deux orbites ouvertes correspondant \`a $\kappa < 0$ et $\kappa > 0$ (pour les d\'etails
voir \cite{hi,hi3}). Quant aux deux autres conditions
(\ref{r2}) et (\ref{c2}), ce sont clairement des conditions \og ouvertes \fg.
\end{proof}

Maintenant si $M$ est une vari\'et\'e presque hermitienne, $V = T_x
M$, $\psi$ est d\'etermin\'e
par le 1-jet de la forme de K\"ahler au point $x$ (en fait seulement par $\omega_x$,
$d\omega_x$) de la fa\c con suivante. Soit $\theta_{\omega}$ la 1-forme telle que
\[ d\omega \wedge \omega = \theta_{\omega} \wedge \omega \wedge \omega
\]
En effet
\beq
\LLambda^1 &     \to & \LLambda^5 \\
\eta      & \mapsto & \eta \wedge \omega \wedge \omega
\eeq
est un isomorphisme lorsque $\omega$ est non d\'eg\'en\'er\'ee. Alors on
pose
\[ \psi_{\omega} = d\omega - \theta_{\omega} \wedge \omega \]
Les \'equations rassembl\'ees dans l'\'enonc\'e de la proposition \ref{open}
portent maintenant sur le 1-jet de $\omega$ : quel que soit $x \in M$ on
appelle $U_x \subset
(\J^1 \LLambda^2)_x$, l'ouvert des 1-jets $j$ tels que
$\omega$, $\psi_{\omega}$ satisfont les conditions
de la proposition \ref{open}.

\begin{theo}
Soit $M$ une vari\'et\'e de dimension 6, $x \in M$. A tout 1-jet $j
\in U_x$, on peut associer localement, sur un voisinage de
$x$, une structure $SU(3)$ telle que la structure $U(3)$ sous-jacente
est de type $W_1 \oplus W_2 \oplus W_4$.
\end{theo}

\begin{proof}
Soit $\omega$ une 2-forme diff\'erentielle dont le 1-jet en $x$, $(j^1
\omega)_x = j$ appartient \`a $U_x$. En un point $y$ suffisamment proche de $x$,
l'op\'eration qui \`a une forme associe son 1-jet \'etant continue,
$(j^1 \omega)_y$ reste dans $U_y$. Il existe un voisinage $N$ de $x$
tel que $(\omega,\psi_\omega)$ v\'erifient (\ref{r1}), (\ref{r2}), (\ref{c2}) en tout
point de $N$. Elles v\'erifient aussi (\ref{c1}) par
d\'efinition de $\psi_{\omega}$, $\theta_{\omega}$ :
\[ \psi_{\omega} \wedge \omega = d\omega \wedge \omega -
\theta_{\omega} \wedge \omega \wedge \omega = 0
\]
Enfin pour que (\ref{c3}) soit v\'erifi\'ee il faut d\'efinir
la structure $SU(3)$ plut\^ot par $(\omega,\psi'=f\psi_{\omega})$ o\`u
\[ f = \sqrt{\frac{2\mu(\omega)}{3\mu(\psi)}} \]
Les \'equations (\ref{r1}) et (\ref{r2}) assurent que $f$ est bien
d\'efinie et que $\omega$, $\psi'$ continuent \`a v\'erifier les
conditions, en particulier $\psi'$ ne s'annule pas. D\`es lors, elles
d\'efinissent une structure $SU(3)$ au voisinage de $x$. En effet, $\mu(f\psi) = f^2 \mu(\psi)$, pour toute fonction $f$, implique $\mu(\psi')= \frac{2}{3}\mu(\omega)$.

Maintenant la vari\'et\'e $(N,g,J)$ est automatiquement de type $W_1
\oplus W_2 \oplus W_4$ par construction :
\[ d\omega = f^{-1}\psi' + \theta_{\omega} \wedge \omega \]
n'a pas de composante dans $\leftr \llambda^{2,1}_0 \rightr$.
\end{proof}

De plus, vue la libert\'e dans le choix de $\omega$, on peut toujours demander que
$d\theta_{\omega} \neq 0$.
\begin{coro}
Il existe des vari\'et\'es presque hermitiennes conformes de type $W_1^c
\oplus W_2^c$ non ferm\'ees.
\label{W124>W12}
\end{coro}

\section*{Conclusion}

La th\'eorie des twisteurs rencontre les vari\'et\'es NK \`a deux
endroits au moins.

Premi\`erement Hitchin \cite{hi2} a d\'emontr\'e que les seuls espaces
de twisteurs au dessus de vari\'et\'es de
dimension 4, k\"ahleriens sont $\CM P^3$ au dessus de $S^4$ et
$\mathbb F^3$, l'espace des drapeaux
de $\CM^3$, au dessus de $\CM P^2$, c'est-\`a-dire 2 parmi les 4 seules
vari\'et\'es homog\`enes de dimension 6 admettant une
structure SNK (voir \cite{bu}). On obtient cette derni\`ere \`a partir
de la structure k\"ahlerienne en
faisant une homoth\'etie et en changeant le signe de la structure presque complexe le long de la fibre, isomorphe \`a $\CM P^1$. De plus, cette m\'ethode est g\'en\'erale pour obtenir des
vari\'et\'es SNK \`a partir de submersions riemanniennes dont
l'espace total est k\"ahlerien, comme a d\'emontr\'e Nagy
dans \cite{na}. Il prouve ainsi l'existence d'une structure SNK sur
l'espace de twisteurs d'une vari\'et\'e K\"ahler-quaternionique.

Deuxi\`emement, par le th\'eor\`eme \ref{twisteurs}, les seules
vari\'et\'es presque hermitiennes conformes de dimension 6 admettant un espace de
twisteurs r\'eduit complexe correspondent aux deux types de
vari\'et\'es NK de dimension 6 : les vari\'et\'es
k\"ahleriennes -- pourvu qu'elles soient Bochner-plates --, et les
vari\'et\'e SNK -- en fait seulement $S^6$ munie de sa structure conforme
standard et de la structure presque complexe compatible issue des octonions. Au
regard de la raison (traduire des probl\`emes sur $M$ dans des
propri\'et\'es des structures holomorphes d'objets associ\'es \`a $Z$, voir
\cite{at}) qui fait s'int\'eresser aux espaces de
twisteurs complexes, ce r\'esultat est juste pr\'eliminaire.

Par exemple, en ce
qui concerne $S^6$, \'etant conform\'ement plate, elle
admet aussi un espace de twisteurs classique. Il s'agit de l'hypersurface
quadrique de dimension 6 complexe $\mathcal Q_+$ (voir \cite{sl}). Par cons\'equent
l'espace de
twisteurs r\'eduit $Z$ est une sous-vari\'et\'e de cette derni\`ere, qu'il
faudra pr\'ecisement d\'ecrire.

En outre, on fait place \`a une observation de
B\'erard-Bergery, Ochiai \cite{ber}. D'une importance cruciale dans la
th\'eorie des twisteurs en dimension 4 est la correspondance
\'etablie en \cite{at} entre les fibr\'es
holomorphes de l'espace de twisteurs, holomorphiquement triviaux sur
chaque fibre, et certains fibr\'es appel\'es
\emph{auto-duaux} sur la base. Cela appelle une
g\'en\'eralisation si possible, en suivant Slupinski \cite{sl} dans le
cas riemannien.

Enfin, une particularit\'e de $S^6$ parmi les
vari\'et\'es SNK de dimension 6 est qu'elle admet plusieurs structures
presque complexes $J$ compatibles avec une m\'etrique donn\'ee
$g$ telles que
$(S^6,g,J)$ est NK. En lien aussi avec la th\'eorie des spineurs de Killing,
on peut esp\'erer donner une expression naturelle de ce fait en
termes de sections de l'espace de twisteurs.

\labelsep .5cm

\end{document}